\documentclass{conm-p-l}
\usepackage{amscd,amssymb,latexsym,subfigure,hyperref}
\usepackage[graph,frame,poly,arc]{xy}  

\theoremstyle{plain}
\newtheorem{theorem}[subsection]{Theorem}

\newtheorem{prop}[subsection]{Proposition}
\newtheorem{cor}[subsection]{Corollary}

\theoremstyle{definition}
\newtheorem{remark}[subsection]{Remark}

\newtheorem{example}[subsection]{Example}
\newtheorem{conj}[subsection]{Conjecture}
\newtheorem{question}[subsection]{Question}
\newtheorem*{note}{Note}

\numberwithin{equation}{section}

\newenvironment{alphenum}%
{%
  \renewcommand{\theenumi}{\alph{enumi}}%
  \renewcommand{\labelenumi}{(\theenumi)}%
  \begin{enumerate}
}
{%
  \end{enumerate}%
}

{%
  \renewcommand{\theenumi}{\roman{enumi}}%
  \renewcommand{\labelenumi}{(\theenumi)}%
  \begin{enumerate}
}
{%
  \end{enumerate}%
}

{%
  \renewcommand{\theenumi}{\Alph{enumi}}%
  \renewcommand{\labelenumi}{(\theenumi)}%
  \begin{enumerate}
}%
{%
  \end{enumerate}%
}

\newenvironment{laundry}
{\begin{list}{$\bullet$}
{
\setlength{\leftmargin}{30pt}
\setlength{\labelsep}{6.5pt}
\setlength{\labelwidth}{4.5pt}
\setlength{\parsep}{1pt}
\setlength{\topsep}{0.3\baselineskip}
\setlength{\itemsep}{0.25\baselineskip}}}
{\end{list}
}

\renewcommand{\arraystretch}{1.1}
\newcommand{\abs}[1]{\lvert #1\rvert}
\newcommand{\bo}[1]{\mathbf{#1}}

\newcommand{\DS}{\displaystyle }
\newcommand{\triv}{\mathbf{1}}
\newcommand{\surj}{\twoheadrightarrow}
\newcommand{\inj}{\hookrightarrow}
\newcommand{\nor}{\triangleleft}

\newcommand{\A}{\mathcal{A}}
\newcommand{\B}{\mathcal{B}}

\newcommand{\LL}{\mathcal{L}}

\newcommand{\RR}{R}

\newcommand{\C}{\mathbb{C}}
\newcommand{\F}{\mathbb{F}}
\newcommand{\Z}{\mathbb{Z}}

\newcommand{\R}{\mathbb{R}}
\newcommand{\N}{\mathbb{N}}
\newcommand{\CP}{\mathbb{CP}}

\newcommand{\K}{\mathbb{K}}

\renewcommand{\SS}{\mathbb{S}}

\newcommand{\PP}{{\mathsf{P}}}
\newcommand{\I}{{\mathfrak{I}}}
\newcommand{\II}{{\check{I}}}

\newcommand{\G}{\Gamma}
\renewcommand{\L}{{\Lambda}}
\renewcommand{\a}{{\alpha}}
\renewcommand{\b}{{\beta}}

\renewcommand{\l}{\lambda}
\newcommand{\s}{\sigma}

\DeclareMathOperator{\rank}{rank}
\DeclareMathOperator{\corank}{corank}
\DeclareMathOperator{\coker}{coker}
\DeclareMathOperator{\id}{id}

\DeclareMathOperator{\Aut}{Aut}
\DeclareMathOperator{\Hom}{Hom}

\DeclareMathOperator{\Epi}{Epi}

\DeclareMathOperator{\GL}{GL}

\DeclareMathOperator{\ab}{ab}

\def\Re{\operatorname{Re}}

\DeclareMathOperator{\Tors}{Tors}

\DeclareMathOperator{\depth}{depth}

\DeclareMathOperator{\ii}{i}
\DeclareMathOperator{\gr}{gr}

\DeclareMathOperator{\ch}{char}
\DeclareMathOperator{\ord}{ord}

\DeclareMathOperator{\Hilb}{Hilb}
\DeclareMathOperator{\dep}{depth}
\DeclareMathOperator{\TC}{TC}
\DeclareMathOperator{\ann}{ann}

\renewcommand{\b}[1]{\mathbf{#1}}
\newcommand{\lin}[1]{{#1}^{\rm lin}}
\newcommand{\ov}[1]{\overline{#1}}
\newcommand{\vsml}[1]{\text{\scriptsize{$#1$}}}
\newcommand{\sml}[1]{\text{\footnotesize{$#1$}}}
\newcommand{\sm}[1]{\text{\small{$#1$}}}

\begin{document}

\title[Fundamental groups of line arrangements]%
{Fundamental groups of line arrangements:\\ Enumerative aspects}

\author{Alexander~I.~Suciu}
\address{Department of Mathematics,
Northeastern University,
Boston, MA 02115}
\email{\href{mailto:alexsuciu@neu.edu}{alexsuciu@neu.edu}}
\urladdr{\href{http://www.math.neu.edu/~suciu/}%
{http://www.math.neu.edu/\~{}suciu}}




\begin{abstract}
This is a survey of some recent developments in the study of 
complements of line arrangements in the complex plane.   
We investigate the fundamental groups and finite covers 
of those complements, focusing on homological and enumerative 
aspects.  

The unifying framework for this study is the stratification 
of the character variety of the fundamental group, $G$, by the jumping 
loci for cohomology with coefficients in rank~$1$ local systems.   
Counting certain torsion points on these ``characteristic" varieties 
yields information about the homology of branched and unbranched 
covers of the complement, as well as on the number of low-index
subgroups of its fundamental group.  

We conclude with two
conjectures, expressing the lower central series quotients of
$G/G''$ (and, in some cases, $G$ itself) in terms of the closely 
related ``resonance" varieties.  We illustrate the discussion with 
a number of detailed examples, some of which reveal new phenomena. 
\end{abstract}

\maketitle

\tableofcontents

\section{Introduction}
\label{sec:intro}

In the introduction to \cite{Hz}, Hirzebruch wrote:  
``The topology of the complement of an arrangement 
of lines in the projective 
plane is very interesting, the investigation of the fundamental 
group of the complement very difficult."  Much progress has 
occurred since that assessment was made in 1983.  The fundamental 
groups of complements of line arrangements are still difficult 
to study, but enough light has been shed on their structure, 
that once seemingly intractable problems can now be attacked 
in earnest.  This paper is meant as an introduction to some 
recent developments, and as an invitation for further investigation.   
We take a fresh look at several topics studied in the past 
two decades, from the point of view of a unified framework.  
Though most of the material is expository, we provide new 
examples and applications, which in turn raise several 
questions and conjectures. 

\subsection{Hyperplane arrangements}
\label{sec:arr}

In its simplest manifestation, an arrangement is merely a finite collection 
of lines in the real plane.  These lines cut the plane into pieces, 
and understanding the topology of the complement amounts to counting 
those pieces.  In the case of lines in the complex plane (or, for that 
matter, hyperplanes in complex $\ell$-space), the complement is connected, 
and its topology (as reflected, for example, in its fundamental group) 
is much more interesting.  

An important example is the braid arrangement of diagonal hyperplanes 
in $\C^{\ell}$.  In that case, loops in the complement can be viewed 
as (pure) braids on $\ell$ strings, and the fundamental group can be 
identified with the pure braid group $P_{\ell}$.  For an arbitrary 
hyperplane arrangement, $\A=\{H_1, \dots , H_n\}$, 
with complement $X(\A)=\C^{\ell}\setminus \bigcup_{i=1}^{n} H_i$, the 
identification of the fundamental group, $G(\A)=\pi_1(X(\A))$, 
is more complicated, but it can be done algorithmically, 
using the theory of braids (see \cite{OT, CSbm}, and the 
references therein).  This theory,
in turn, is intimately connected with the theory of knots and links 
in $3$-space, with its wealth of algebraic and combinatorial invariants, 
and its varied applications to biology, chemistry, and physics. 
A revealing example where developments in arrangement 
theory have influenced knot theory is Falk and Randell's \cite{FR88} 
proof of the residual nilpotency of the pure braid group, 
a fact that has been 
put to good use in the study of Vassiliev invariants.  We refer to 
the excellent surveys \cite{FR86, FR00} for a more 
complete treatment of the homotopy theory of arrangements. 

A more direct link to physics is provided by the deep  
connections between arrangement theory and hypergeometric functions. 
Work by Schechtman-Varchenko \cite{SV} and many others has  
profound implications in the study of Knizhnik-Zamolod\-chikov equations 
in conformal field theory. 
We refer to the recent monograph by Orlik and Terao \cite{OT00} for 
a comprehensive account of this fascinating subject. 

Hyperplane arrangements, 
and the closely related configuration spaces, are used 
in numerous areas, including robotics, graphics, molecular 
biology, computer vision, and databases for representing 
the space of all possible states of a system characterized 
by many degrees of freedom.  Understanding the topology 
of complements of subspace arrangements and configuration spaces 
is important in robot motion planning (finding a collision-free 
motion between two placements of a given robot among a set of 
obstacles, see \cite{Gh}), and in multi-dimensional billiards 
(describing periodic trajectories of a mass-point in a domain 
in Euclidean space with reflecting boundary, see \cite{FT}).   

\subsection{Characteristic and resonance varieties}
\label{sec:cvs}
The unifying framework for the study of fundamental groups of 
hyperplane arrangements is provided by their characteristic 
and resonance varieties.  The origins of 
those varieties can be traced back to the work of Novikov \cite{Nov} 
in the mid 1980's 
on the cohomology of smooth manifolds with coefficients in local systems, 
and its relation to Morse theory for $1$-forms.  Novikov showed 
that the Betti numbers remain constant almost everywhere on the 
representation space of the fundamental group of the manifold, 
but they increase on the union of countably many algebraic 
submanifolds.

The simplest situation is when the fundamental group 
of the manifold, $G=\pi_1(X)$,  
has abelianization $H_1(G)\cong\Z^n$, and the rank 
of the local coefficient system is $1$. 
Then, the representation variety is the algebraic torus  
$\Hom(G,\C^*)\cong {\C^*}^n$, and the jumping loci for 
$1$-dimensional cohomology are the {\em characteristic varieties} 
of $G$, 
\begin{equation*} 
\label{eq:cvdef}
V_d(G,\C)=\{\bo{t}\in {\C^{*}}^n \mid 
\dim_{\C} H^1(G,\C_{\mathbf{t}}) \ge d\},\quad 1\le d \le n.
\end{equation*}
Work of Libgober \cite{Li1}, Hironaka \cite{Hi1, Hi2}, 
and Sakuma \cite{Sa} clarified the connection between 
the characteristic varieties, the Alexander matrix, 
and the Betti numbers of finite abelian covers.\footnote{%
In fact, one can argue that the origins of characteristic varieties 
(so named in \cite{Li2}) 
go back to the work of James Alexander, in the early 1920's.  
Indeed, for a knot in $\SS^3$, the characteristic varieties 
are nothing but the roots of the Alexander polynomial; 
if the knot is fibered, those are the eigenvalues of 
the algebraic monodromy, whence the name for those varieties  
(also known as Alexander varieties \cite{Hi2}).}

A breakthrough occurred in 1997, when Arapura~\cite{Ar} showed that 
the characteristic varieties of $G=\pi_1(X)$ are unions of 
(possibly torsion-translated) subtori of ${\C^*}^{n}$, 
provided $X$ is the complement of a normal-crossing divisor 
in a compact K\"{a}hler manifold with vanishing first Betti 
number.\footnote{%
Arapura's theorem works for higher-rank representations of $\pi_1(X)$, 
as well as for the jumping loci of higher-dimensional cohomology 
groups of $X$.}
In particular, this result applies to the characteristic varieties 
of an arrangement group, $G=G(\A)$. 
Also in 1997, Falk \cite{Fa} defined the 
{\em resonance varieties} of the Orlik-Solomon algebra 
$H^*(X(\A),\C)$.   There ensued a flurry of activity, 
showing (in \cite{CScv, Li2, CO, Li3}) that these varieties are the 
tangent cones at the origin to the characteristic varieties.  
The combinatorial description of the resonance varieties of arrangements, 
started by Falk, was completed by Libgober and Yuzvinsky \cite{LY}. 

\subsection{Betti numbers of finite covers}
\label{sec:betticov}
One of the reasons for studying the characteristic 
varieties of the group $G=\pi_1(X)$ is the precise information 
they give about the (rational) homology of covering spaces of $X$.  
Indeed, suppose $X'\to X$ is a regular cover,  
with finite, abelian Galois group.  
Results of Libgober~\cite{Li1}, Sakuma~\cite{Sa}, and Hironaka~\cite{Hi3}  
show how to compute $b_1(X')$ by counting torsion points 
of a certain order on the algebraic torus $\Hom(G,\C^*)$, 
according to their depth in the filtration by the varieties $V_d(G,\C)$.  
(Explicit formulas of this sort are given in 
Corollary~\ref{cor:torsdepth} and Theorem~\ref{thm:reflibsak}.) 

Of particular interest is the case of {\em congruence covers}, 
$X_N$, defined by the canonical projection $G\to H_1(G,\Z_N)$. 
As noted by Sarnak and Adams \cite{SarA}, the sequence of Betti numbers 
$\{b_1(X_N)\}_{N\in \N}$ is {\em polynomially periodic}, i.e., 
there exists polynomials $P_1(x),\dots ,P_{T}(x)$, such that 
$b_1(X_N)=P_i(N)$, whenever $N\equiv i\,\bmod T$.
(This follows from Sakuma's formula, by means of a 
deep result of Laurent on torsion points on varieties). 
Here we give examples of arrangements $\A$ for which 
the minimal period, $T$, is greater than $1$ (in fact, $T=2$).  
This non-trivial periodicity in the Betti numbers of congruence covers 
of $X(\A)$ can be traced back to the presence of 
isolated torsion points in $V_2(G(\A),\C)$.  

Similar techniques apply to the study of {\em Hirzebruch covering 
surfaces\footnote{%
Not to be confused with the {\em Hirzebruch surfaces}, defined 
30 years previously.}}, $M_N(\A)$, defined in \cite{Hz}, as follows.
Let $\A$ be a central arrangement of planes 
in $\C^3$, and let $\ov{\A}$ be the corresponding arrangement 
of lines in $\CP^2$.  Then $M_N(\A)$ is the minimal desingularization 
of the branched congruence cover, $\widehat{X}_N(\ov{\A})$, of the 
projective plane, ramified along $\ov{\A}$.  
The polynomial periodicity of the sequence 
$\{b_1(M_N(\A)\}_{N\in \N}$ was established by Hironaka \cite{Hi0} 
and Sakuma \cite{Sa}.  We give examples of arrangements 
for which the minimal period is $T=4$.   The non-trivial periodicity 
in the Betti numbers of Hirzebruch covering surfaces is due 
to the presence of torsion-translated, positive-dimensional 
components in $V_1(G(\A),\C)$.  

\subsection{Characteristic varieties over finite fields}
\label{sec:finitecv}
The usefulness of characteristic varieties 
in the study of finite abelian covers $X'\to X$ 
is limited by the fact that they only give the rank 
of $H_1(X')$, but not its torsion coefficients.  
This limitation can be overcome (at least partially) 
by considering the characteristic varieties of $G=\pi_1(X)$  
over arbitrary fields $\K$.  These varieties were first 
considered by Matei in his thesis~\cite{Mt}; the analogous resonance 
varieties over finite fields were first studied in \cite{MS2}. 
From the stratification of $\Hom(G,\K)$ by the varieties 
$V_d(G,\K)$, one can derive information about the 
torsion in $H_1(X')$, away from the prime $q=\ch\K$.  
Such a method (which recovers and extends results of 
Libgober, Sakuma, and Hironaka) was recently developed in \cite{MS3}.   
We sketch this method here, and mention some of 
its applications, which show why it is important to look  
beyond Betti numbers. 

In \cite{Fx70}, Fox showed how to compute the number of 
metacyclic representations of link groups, by means of 
his free differential calculus.  In \cite{MS3}, we generalize 
Fox's method, and interpret it in terms of characteristic 
varieties over certain Galois fields $\K$.  This provides 
a way to count homomorphisms from a finitely-presented 
group $G$ onto certain finite metabelian groups, such as 
the symmetric group $S_3$ and the alternating group $A_4$. 
In turn, this count gives information about the number 
of finite-index subgroups of $G$, 
at least for low index.  All these enumeration problems 
reduce to counting torsion points on the characteristic 
varieties $V_d(G,\K)$.  Such a count depends in subtle ways 
on the field $\K$, and is certainly affected by the presence 
of translated subtori in $V_d(G,\C)$.  

\subsection{Lower central series}
\label{sec:lcs}
Probably the most studied numerical invariants of 
an arrangement group $G$ are the ranks of its lower central 
series quotients, $\phi_k(G)=\rank \gamma_k G/\gamma_{k+1} G$. 
The impetus came from the work of Kohno~\cite{Ko}, who used 
rational homotopy theory to compute $\phi_k(P_{\ell})$, 
the LCS ranks of the braid arrangement group.   
Falk and Randell \cite{FR1}, using more direct methods, 
established the celebrated {\em LCS formula} for the 
broader class of {\em fiber-type} arrangements, expressing 
$\phi_k(G(\A))$ in terms of the exponents of the arrangement $\A$. 
Recent work by a number of authors, \cite{SY, CScc, JP, PY, Pe},  
has further clarified the meaning of this formula, and greatly 
expanded its range of applicability.  

Another direction was started in \cite{CSpn, CSai}, with the 
study of the ranks of Chen groups of arrangements, 
$\theta_k(G)=\phi_k(G/G'')$. The $\theta$-invariants 
can provide stronger information than the 
$\phi$-invariants, distinguishing, in some cases, 
groups of fiber-type arrangements from the corresponding 
direct products of free groups.  

A number of experiments (some of which we record here) has 
revealed tantalizing parallels to other, seemingly unrelated 
objects associated to the group $G$---namely, its resonance 
varieties. As a result, we formulate two conjectures, expressing the 
ranks of the LCS quotients of $G/G''$ (and, in some cases, 
$G$ itself) in terms of the dimensions of the components of 
the resonance variety $\RR_1(G,\C)$.

\subsection{Organization of the paper}
\label{sec:org}
Sections~\ref{sec:pi1arr} and \ref{sec:rvcv} contain basic 
material on complements of complex hyperplane arrangements, 
their fundamental groups, Alexander matrices, and resonance 
and characteristic varieties.  Much of the exposition is 
based on joint work with D.~Cohen~\cite{CSbm, CSai, CScv} 
and D.~Matei~\cite{MS1, MS2, MS3}. 

In Section~\ref{sec:fincov}, we discuss the connection 
between the characteristic varieties of a group  
and the homology of its finite-index subgroups.  
We use this connection in Section \ref{sec:counting} to count 
low-index subgroups.  These two sections are a partial summary  
of a recent paper with D.~Matei~\cite{MS3}. 

In Sections~\ref{sec:cong} and \ref{sec:hirzebruch}, 
we study congruence covers and Hirzebruch covering surfaces, 
and the polynomial periodicity of their Betti numbers.

In Section~\ref{sec:lcschen}, we discuss the LCS quotients  
and Chen groups of an arrangement group, and their 
possible connection to the resonance varieties. 

Section~\ref{sec:examples} contains about a dozen examples 
worked out in detail.  A guide to the examples is provided in 
Section~\ref{sec:glossary}.  

\subsection*{Acknowledgments}  
I wish to thank Dan Cohen, Mike Falk, Anatoly Libgober, 
Daniel Matei, \c{S}tefan Papadima, Dick Randell, and 
Sergey Yuzvinsky for useful discussions and comments 
while some of this material was being developed.  
Especially, Dan Cohen, who attracted me to the subject, 
and Daniel Matei, together with whom a good deal 
of the work presented here was done. 

I also wish to thank Emma Previato, whose encouragement 
and enthusiasm convinced me to undertake this project, 
and the referee, for pertinent remarks.
 
A preliminary version of the paper was presented in the 
Special Session on Enumerative Geometry in Physics, 
at the AMS Spring Eastern Section Meeting, held in 
Lowell, MA in April 2000.  A fuller version was  
presented as an Invited Address at the AMS Fall Eastern 
Section Meeting, held in New York, NY in November 2000.


\section{Fundamental groups and Alexander matrices}
\label{sec:pi1arr}

\subsection{Arrangements of hyperplanes}
\label{sec:hyparr}

A {\em (complex) hyperplane arrangement} is a finite
set, $\A$, of codimension-$1$ affine subspaces in a 
finite-dimensional complex vector space, $\C^{\ell}$.  
The main combinatorial object associated to $\A$ 
is its {\em intersection lattice}, 
$\LL(\A)=\left\{\emptyset \ne \bigcap_{H\in \B}\, H\left| 
\B\subseteq\A\right\}\right.$. 
This is a ranked poset, consisting of all non-empty intersections
of $\A$, ordered by reverse inclusion, and with rank function given 
by codimension.  The main topological object associated to $\A$ is 
its {\em complement}, $X(\A)=\C^{\ell}\setminus\bigcup_{H\in \A}H$. 
This is an open, $\ell$-dimensional complex manifold, whose topological 
invariants are intimately connected to the combinatorics of the arrangement. 
A basic reference for the subject is the book by Orlik and Terao~\cite{OT}.

\begin{example}
\label{ex:braidarr}
Probably the best-known hyperplane arrangement 
is the {\em braid arrangement} in $\C^{\ell}$.  
This is the complexification of the 
$\operatorname{A_{\ell-1}}$-reflection arrangement,   
$\B_{\ell}=\{\ker (z_i-z_j)\}_{1\le i<j\le \ell}$, with   
lattice $\LL(\B_{\ell})=\Pi(\ell)$, the partition 
lattice of $[\ell]:=\{1,\dots,\ell\}$, complement 
$X(\B_{\ell})=F(\ell,\C)$, the configuration space 
of $\ell$ ordered points in $\C$, and group 
$\pi_1(X(\B_{\ell}))=P_{\ell}$, the pure braid group 
on $\ell$ strings.
\end{example}

Let $\A=\{H_1,\dots , H_n\}$ be a {\em central} arrangement in $\C^{\ell}$ 
(that is, $\b{0}\in \bigcap_{i} H_i$).           
A defining polynomial for $\A$ may be written as $Q_{\A}=f_1\cdots f_n$, 
where $f_i$ are distinct linear forms, with $H_i=\ker f_i$. Choose coordinates
$(z_1,\dots,z_{\ell})$  in $\C^{\ell}$ so that $H_n =\{z_{\ell}=0\}$.  
The corresponding {\em decone} of $\A$ 
is the affine arrangement $\A^*=\mathbf{d}\A$ in
$\C^{\ell-1}$,  with defining polynomial 
$Q_{\A^*} = Q_{\A}(z_1,\dots,z_{\ell-1},1)$.  
Reversing the procedure yields $\A=\mathbf{c}\A^*$,  
the {\em cone} of $\A^*$.  
The respective complements are related as follows: 
$X(\A) \cong X(\A^*) \times\C^*$, where $\C^*=\C\setminus \{ 0 \}$.   

\subsection{Arrangements of lines}
\label{subsec:lines}
We will be interested here in the low-dimen\-sional topology of 
arrangements---mainly in questions related to the fundamental group 
of the complement, and the cohomology ring up to degree~$2$.   
For that, we may restrict our attention to (affine) line arrangements 
in $\C^2$, or, by coning, to (central) plane arrangements in $\C^3$.  
Indeed, if $\A$ is a hyperplane arrangement, with complement $X$, 
let $\A'$ be a generic two-dimensional section of $\A$, with complement $X'$.  
Then, by the Lefschetz-type theorem of Hamm and L\^{e}~\cite{HL}, 
the inclusion $j:X'\to X$ induces an isomorphism 
$j_*:\pi_1(X')\to\pi_1(X)$ and a monomorphism   
$j^*:H^{2}(X)\to H^{2}(X')$.  By the Brieskorn-Orlik-Solomon 
theorem (cf.~\cite{OT}), the map $j^*$ is, in fact, an isomorphism, 
and thus $H^{\le 2}(X)\cong H^{\le 2}(X')$.  
Furthermore, the Hurewicz homomorphism $\pi_2(X)\to H_2(X)$ 
is the zero map, see Randell~\cite{Ra2}.  
Thus, by a theorem of Hopf, $H^{\le 2}(\pi_1(X))\cong H^{\le 2}(X)$,    
see \cite{MS2}.  

So let $\A =\{l_{1},\dots ,l_{n}\}$ be an arrangement of $n$ affine lines
in $\C^2$.  Let $v_1, \dots ,v_s$ be the intersection points of the lines. 
If $v_q = l_{i_1}\cap \dots \cap l_{i_r}$, set 
$I_q=\{i_1,\dots,i_{r  }\}$.  
Let $\LL_1(\A)=[n]$ and $\LL_2(\A)=\{I_{1},\dots,I_{s}\}$. 
The combinatorics of the arrangement is encoded in its intersection 
poset, $\LL(\A)=\{\LL_1(\A),\LL_2(\A)\}$, which keeps track of the
incidence relations between the points and the lines of the
arrangement. 

From $\LL(\A)$, one can extract some simple, but 
very useful numerical information.  
First note that $n=\abs{\LL_1(\A)}$, $s=\abs{\LL_2(\A)}$.  
Let 
\begin{equation} 
\label{eq:multiplicities}
m_r=\#\{ I\in \LL_2(\A) \mid \abs{I}=r\}
\end{equation} 
be the number of intersection points of multiplicity $r$. 
Note that $\sum_{r} m_r=s$.  
The homology groups of the complement, $X=X(\A)$, 
are easily computable in terms of these data  
(see~\cite{OT}): 
\begin{equation} 
\label{eq:homology}
H_1(X)=\Z^n,\quad H_2(X)=\Z^{\sum_{r} m_r(r - 1)}, \quad 
H_i(X)=0,\ \text{for $i>2$}.
\end{equation}

\subsection{Fundamental groups}
\label{subsec:pi1arr}
We now describe a procedure for finding a finite 
presentation for the fundamental group of the complement 
of an arrangement, 
$G(\A)=\pi_1(X(\A))$.  The information encoded in the 
intersection lattice is not {\it a priori} enough; we need 
slightly more information, encoded in the {\em braid monodromy} 
of the arrangement.   As shown in \cite{CSbm}, the resulting 
presentation is Tietze-I equivalent to the presentations of 
Randell \cite{Ra} and Arvola \cite{Ar}.

Let $B_n$ be the Artin braid group on $n$ strings, with 
generators $\sigma_1,\dots,\sigma_{n-1}$.   
Let $P_{n}=\ker(B_n\to S_n)$ be the subgroup of braids
with trivial permutation of the strings, with generators 
$A_{ij}=\sigma _{j-1}\cdots \sigma _{i+1}
\sigma _{i}^{2}\sigma _{i+1}^{-1}\cdots \sigma _{j-1}^{-1}$ 
($1\le i<j\le n$).    
The pure braid group acts on the free group 
$F_{n}=\langle x_{1},\dots ,x_{n}\rangle $ 
by restricting the Artin representation $B_n\inj \Aut(F_n)$ to $P_n$.  
Explicitly:
\begin{equation} 
\label{eq:artin}
A_{ij} (x_r) =
\begin{cases}
x_r&\text{if $r<i$ or $r>j$,}\\
x_ix_rx_i^{-1}&\text{if $r=j$,}\\
x_ix_jx_rx_j^{-1}x_i^{-1}&\text{if $r=i$,}\\ \relax
[x_i,x_j]x_r[x_i,x_j]^{-1}&\text{if $i<r<j$,}
\end{cases}
\end{equation}
see Birman~\cite{Bi}.
For an increasingly ordered set $I=\{i_{1},\dots ,i_{r}\}\subset [n]$, 
let $\II=\{i_1,\dots ,i_{r-1}\}$, and let 
\begin{equation}
\label{eq:AI}
A_{I} = A_{i_{1}i_{2}} A_{i_{1}i_{3}}A_{i_{2}i_{3}} 
A_{i_{1}i_{4}}A_{i_{2}i_{4}}A_{i_{3}i_{4}} \cdots 
A_{i_{1}i_{r}}\cdots A_{i_{r-1}i_{r}} 
\end{equation}
be the braid in $P_n$ which 
performs a full twist on the strands corresponding to $I$, leaving 
the other strands fixed. 

Now choose a generic linear projection $p:\C^2\to \C$, and a 
point $y_0\in \C$, with $\Re(y_0)\gg 0$.  
Let $l_0$ be the line $p^{-1}(y_0)$. 
Label the lines $l_i$ and the vertices $v_q$ so that 
$\Re(p(v_1))>\cdots>\Re(p(v_s))$ and 
$\Re(l_0\cap l_1)<\cdots <\Re(l_0\cap l_n)$. 
Pick a basepoint $\tilde{y}_0\in \C^2$ so that 
$\tilde{y}_0\in l_0$, but $\tilde{y}_0\notin l_i$. 
Let $x_i$ be the meridian loops to $l_i$, 
based at $\tilde{y}_0$, and oriented compatibly with 
the complex orientations of $l_i$ and $\C^2$. 
There are then pure braids 
\begin{equation}
\label{eq:braidmonogens}
\a_1=A_{I_{1}}^{\delta_{1}},\dots ,\a_s=A_{I_{s}}^{\delta_{s}}, 
\end{equation}
where $a^b:=b^{-1}ab$, such that $G=G(\A)$ has braid monodromy 
presentation
\begin{equation}
\label{eq:bmpres}
G=G(\a_1,\dots,\a_s):=\langle x_1,\dots ,x_n \mid \alpha_{q}(x_i)=x_i
\ \text{for } i \in \II_q \text{ and } q\in [s]
\rangle.  
\end{equation}

The conjugating braids $\delta_{q}$ may be obtained  as follows. 
In the case where $\A$ is the complexification of a real arrangement, 
each vertex set $I_{q}$  gives rise to a partition
$[n] = I'_{q}\cup I_{q} \cup I''_{q}$ into lower, middle, and upper indices.  
Let $J_{q}=\{i \in I''_{q} \mid \min I_{q}<i<\max I_{q}\}$.  
Then $\delta_{q}$ is the subword of the full twist 
$A_{12\dots n}=\prod_{i=2}^n\prod_{j=1}^{i-1} A_{ji}$, 
given by 
\begin{equation}
\label{eq:conjbraid}
\delta_{q} = \prod_{i\in I_{q}} \prod_{j\in J_{q}} A_{ji},
\end{equation}
see \cite{Hi1, CF, CSbm}.   
In the general case, the braids $\delta_q$ can be 
read off a ``braided wiring diagram"  associated to $\A$ and the 
projection $p$, see \cite{CSbm} for further details.

\begin{figure}[ht]
\subfigure{%
\label{fig:4lines-a}%
\begin{minipage}[t]{0.5\textwidth}
\setlength{\unitlength}{17pt}
\begin{picture}(4.5,4.5)(-2.5,-0.5)
\put(0,0){\line(1,1){4}}
\put(-1,2){\line(1,0){6}}
\put(0,4){\line(1,-1){4}}
\put(-0.5,0.7){\line(3,1){5.5}}
\put(4.4,-0.5){\makebox(0,0){$l_1$}}
\put(5.5,1.95){\makebox(0,0){$l_2$}}
\put(5.5,2.7){\makebox(0,0){$l_3$}}
\put(4.4,4.5){\makebox(0,0){$l_4$}}
\put(2,2){\makebox(0.15,0){\circle*{0.14}}}
\put(3.4,2){\makebox(0.17,0){\circle*{0.14}}}
\put(2.4,1.65){\makebox(0.04,0){\circle*{0.14}}}
\put(1.37,1.3){\makebox(0,0){\circle*{0.14}}}
\put(3.4,2.3){\makebox(0,0){\sm{v_1}}}
\put(2.4,1.2){\makebox(0,0){\sm{v_2}}}
\put(2,2.5){\makebox(0,0){\sm{v_3}}}
\put(1,1.6){\makebox(0,0){\sm{v_4}}}
\end{picture}
\end{minipage}
}
\subfigure{%
\label{fig:4lines-b}%
\begin{minipage}[t]{0.3\textwidth}
\setlength{\unitlength}{18pt}
\begin{picture}(4,4)(0,-0.5)
\xygraph{!{0;<10mm,0mm>:<0mm,20mm>::}
[]*D(3){I_1}*-{\bullet}
(-[dr]*U(3){2}*-{\bullet}
,-[drr]*U(3){3}*-{\bullet}
,[r]*D(3){I_2}*-{\bullet}
(
,-[dl]*U(3){1}*-{\bullet}
,-[dr]
)
,[rr]*D(3){I_3}*-{\bullet}
(
,-[dll]
,-[dl]
,-[dr]*U(3){4}*-{\bullet}
)
,[rrr]*D(3){I_4}*-{\bullet}
(
,-[dl]
,-[d]
)
)
}
\end{picture}
\end{minipage}
}
\caption{\textsf{An arrangement of $4$ lines 
in $\C^2$ and its intersection poset}}
\label{fig:toy}
\end{figure}

\begin{example}
\label{ex:toy}
Here is a simple example which illustrates the procedure.  
Let $\A$ be the arrangement of $4$ lines in $\C^2$ defined 
by $Q_{\A}=z(z-y)(z+y)(2z-y+1)$.  
The lines $l_1,l_2,l_3,l_4$, together with  the vertices 
$v_1,v_2,v_3,v_4$ are depicted in Figure~\ref{fig:toy}. 
The intersection lattice is also shown.  
Take $p:\C^2\to \C$ to be the 
projection $p(y,z)=y$, and choose $y_0=2$. 
The corresponding braid monodromy generators are:  
$\a_1=A_{23}$, $\a_2=A_{13}^{A_{23}}$, $\a_3=A_{124}$, 
$\a_4=A_{34}$.  Replacing $\a_2$ by $\a_1\a_2 \a_1^{-1}=A_{13}$ 
gives an equivalent presentation for $G=\pi_1(X(\A))$.  We obtain:
\begin{align*}
G&=G(A_{23}, A_{13}, A_{124}, A_{34})\\
&=\langle x_1,x_2,x_3,x_4 \mid x_1x_2x_4=x_4x_1x_2=x_2x_4x_1, 
[x_1,x_3]=[x_2,x_3]=[x_4,x_3]=1\rangle.
\end{align*}
Consequently, $G\cong F_2\times \Z^2$.
\end{example}

\subsection{Fox calculus}
\label{subsec:fox}
Consider a finitely-presented group $G$, with  
presentation $G=\langle x_1,\dots ,x_{n}\mid r_1,\dots ,r_m\rangle$.    
Let $X_G$ be the $2$-complex modelled on this presentation,  
$\widetilde{X}_G$ its universal cover, and $C_*(\widetilde{X}_G)$  
its augmented cellular chain complex. Picking as
generators for the chain groups the lifts of the cells
of $X_G$, the complex $C_*(\widetilde{X}_G)$ becomes identified with
\begin{equation}
\label{eq:chaincomplex}
(\Z G)^m \xrightarrow{
J_G=
\begin{pmatrix}
\frac{\partial r_i}{\partial x_j}\end{pmatrix}^{\phi}
} (\Z G)^{n}
\xrightarrow{
\vsml{
\begin{pmatrix}
 x_1-1\\ \cdots  \\ x_{n}-1
\end{pmatrix}
}
}\Z G\xrightarrow{\epsilon}\Z\to 0,   
\end{equation}
where $\epsilon$ is the augmentation map,  $\phi:F_n\to G$ 
is the projection associated to the presentation, 
and $\frac{\partial}{\partial x_j}: \Z F_{n}\to \Z F_{n}$ 
are the Fox derivatives, defined by the rules
\begin{equation} 
\label{eq:forder}
\frac{\partial 1}{\partial x_j}=0,\quad
\frac{\partial x_i}{\partial x_j}=\delta_{ij}, \quad
\frac{\partial(uv)}{\partial x_j}=\frac{\partial u}{\partial x_j}\epsilon(v)+
u\frac{\partial v}{\partial x_j}.
\end{equation}

Now assume that $H_1(G)\cong \Z^n$.  
Let $\ab:G\to \Z^n$ be the abelianization map, and set $t_i=\ab(x_i)$. 
This choice of generators for $\Z^n$ identifies the group ring $\Z\Z^n$   
with the ring of Laurent polynomials $\Lambda_n=\Z[t_1^{\pm 1}, \dots,t_n^{\pm 1}]$
The {\em Alexander matrix}	of $G$ is the $m\times n$ matrix $A_G$, with entries in 
$\Lambda_n$, obtained by abelianizing the Fox Jacobian of $G$ (see \cite{Fx1, Bi}): 
\begin{equation} 
\label{eq:alexmat}
A_G=J_G^{\ab}.
\end{equation}  

Let $\psi:\Lambda_n\to
\Z[[\l_1,\dots, \l_n]]$ be the ring homomorphism given by $\psi(t_i)=1-\l_i$ 
and $\psi(t_i^{-1})=\sum_{q\ge 0}\l_{i}^{q}$, and let $\psi^{(h)}$ 
be its homogeneous part of degree $h$. Since $H_1(G)\cong \Z^n$, 
the entries of $A_G$ are in the augmentation ideal 
$\I=(t_1-1,\dots, t_n-1)$, and so $\psi^{(0)}A_G$ is the zero matrix.  
The {\em linearized Alexander matrix} of $G$ is the $m\times n$ matrix 
\begin{equation} 
\label{eq:linalexmat}
\lin{A}_G=\psi^{(1)}A_G,
\end{equation}  
obtained by taking the degree $1$ part of 
the Alexander matrix. Its entries are integral linear forms in 
$\l_1,\dots, \l_n$.  As noted in \cite{MS2}, the coefficient 
matrix of  $\lin{A}_G$ can be interpreted 
as the augmented Fox Hessian of $G$: 
\begin{equation} 
\label{eq:hesse}
(\lin{A}_G)_{i,j}=\sum_{k=1}^{n}\epsilon \Big(
\frac{\partial^2 r_i}{\partial x_k \partial x_j} \Big) \l_k.
\end{equation}

\subsection{Alexander matrices of arrangements}
\label{subsec:alex}
Let $\A$ be an arrangement, with $\abs{\A}=n$, 
$\LL_2(\A)=\{I_1,\dots ,I_s\}$, and braid monodromy generators 
$\a_q=A_{I_q}^{\delta_q}\in P_n$.
The Alexander matrix of $G(\A)=G(\a_1,\dots ,\a_s)$ can be computed directly 
from the Gassner representation of the braid monodromy generators, 
whereas the linearized Alexander matrix can 
be computed directly from the intersection lattice:

\begin{alphenum}
\item 
The {\em Gassner representation} of the pure braid group, 
$\Theta:P_n\to \GL(n,\L_n)$, is defined by 
\begin{equation}
\label{eq:gass}
\Theta(\a)=\begin{pmatrix}
\DS{\frac{\partial \a(x_i)}{\partial x_j}}
\end{pmatrix}^{\ab}.
\end{equation}

\item
The matrix $A_G$ is obtained by stacking the matrices 
$\Theta (\a_{1})-\id, \dots$, 
$\Theta (\a_{s})-\id$, 
and selecting the rows corresponding to $\II_1,\dots ,\II_s$. 

\item
The  matrix $\lin{A}_G$ is obtained by 
stacking the matrices $\lin{A}_{I_1}, \dots, \lin{A}_{I_s}$, where
\begin{equation}
\label{eq:linalex}
(\lin{A}_{I})_{i,j}= \delta_{j,I} \Big(\lambda_i -
\delta_{i,j} \sum_{k\in I} \lambda_k \Big),\quad
\text{for } i\in \II \text{ and } j\in [n].
\end{equation}
and $\delta_{j,I}=1$ if $j\in I$, and $0$ otherwise.
\end{alphenum}

\begin{example}
\label{ex:alextoy}
For the ``near-pencil" of $4$ lines in Example~\ref{ex:toy}, 
the Alexander matrix, and its linearization, are as follows:
\begin{align*}
A &=
\begin{pmatrix}
t_{1}(t_{2}t_{4}-1)
&t_{1}(1-t_{1})
&0
&t_{1}t_{2}(1-t_{1})\\
1-t_{2}
&t_{1}t_{2}(t_{4}-1)+t_{1}-1
&0
&t_{1}t_{2}(1-t_{2})\\
t_{1}(t_{3}-1)
&0
&t_{1}(1-t_{1})
&0\\
0
&t_{2}(t_{3}-1)
&t_{2}(1-t_{2})
&0\\
0
&0
&t_{3}(t_{4}-1)
&t_{3}(1-t_{3})
\end{pmatrix},
\\
\lin{A} &=
\begin{pmatrix}
-\l_2-\l_4 & \l_1 & 0 & \l_1 \\
\l_2  & -\l_1-\l_4 & 0 & \l_2 \\
-\l_3 & 0 & \l_1 & 0 \\
0 & -\l_3 & \l_2 & 0 \\
0 & 0 & -\l_4 & \l_3
\end{pmatrix}.
\end{align*}

\end{example}

\section{Resonance varieties and characteristic varieties}
\label{sec:rvcv}

\subsection{Resonance varieties}
\label{subsec:rv}

Let $G=\langle x_1,\dots ,x_{n}\mid r_1,\dots ,r_m\rangle$ 
be a finitely-presented group. Let $H^*(G,\K)$ be the cohomology 
ring of $G$, with coefficients in a field $\K$ (with cup-product 
denoted by $\cdot$).   

The {\em $d$-th resonance variety} of $G$, over the 
field $\K$, is the set of $\l\in  H^1(G,\K)$ for which 
there exists a subspace $W\in  H^1(G,\K)$ of dimension $d+1$ 
such that $\mu\cdot \lambda =0$, for all $\mu\in W$. 
If $\ch\K\ne 2$, we have $\lambda\cdot\lambda=0$, and thus 
\begin{equation} 
\label{eq:resvar}
\RR_d(G,\K)=\{\l\in H^1(G,\K) \mid 
\dim_{\K} H^1(H^*(G,\K),\cdot \l) \ge d\}, 
\end{equation}
where $(H^*(G,\K),\cdot \l)$ is the chain complex with chains 
the cohomology groups of $G$, and differentials    
equal to multiplication by $\l$.

If $H_1(G)\cong \Z^n$, then $H^1(G,\K)$ may be identified with $\K^n$.  
If, moreover, 
$H_2(G)$ is torsion-free, then $\RR_{d}(G,\K)$ may be identified with 
the $d$-th determinantal variety of the linearized Alexander matrix: 
\begin{equation} 
\label{eq:rlinal}
\RR_{d}(G,\K)=\{\l\in \K^n \mid \rank_{\K} \lin{A}_G(\l) < n-d \},
\end{equation}
where $\lin{A}_G(\l)$ denotes the linearization of $A_G$, 
evaluated at $\l=(\l_1,\dots,\l_n)\in\K^n$, 
see \cite[Theorem~3.9]{MS2}.

The $\K$-resonance varieties form a descending filtration 
$\K^n=\RR_0\supset \RR_1\supset\cdots\supset \RR_{n-1}\supset \RR_n=\{\b{0}\}$. 
The ambient type of each term in the filtration depends only on the
truncated cohomology ring $H^{\le 2}(G)$. More precisely, if 
$H^{\le 2}(G_1) \cong H^{\le 2}(G_2)$, there exists a linear automorphism 
of $\K^n$ taking $\RR_{d}(G_1,\K)$ to $\RR_{d}(G_2,\K)$.

\subsection{Resonance varieties of arrangements}
\label{subsec:rvarr}

The resonance varieties of a complex hyperplane arrangement $\A$ were 
first defined by Falk in~\cite{Fa}, as the cohomology jumping loci of the 
Orlik-Solomon algebra, $A(\A)=H^*(X(\A),\C)$.  In particular, for an arrangement 
$\A$ of $n$ hyperplanes, he defined $\RR_d(\A)$ to be the set of 
$\l\in H^1(X(\A),\C)\cong \C^n$ for which 
$\dim_{\C} H^1(A(\A),\cdot \l) \ge d$.  
As shown in \cite{MS2}, this definition agrees with the one 
in~\eqref{eq:resvar}, i.e, $\RR_d(\A)=R_d(G,\C)$, where $G=G(\A)$.  

The (complex) resonance varieties of arrangements are by now 
very well understood (see the recent surveys by Falk \cite{Fa00} 
and Yuzvinsky \cite{Yu00} for a thorough treatment).  
As noted by Falk~\cite{Fa}, $\RR_1(\A)$ is contained in the hyperplane 
\begin{equation} 
\label{eq:deltan}
\Delta_n=\{\l\in \C^n \mid \sum_{i=1}^n\l_i=0\}.
\end{equation}
As shown in \cite{CScv}, each component of $\RR_1(\A)$ 
is a linear subspace of $\C^n$ (see \cite{Li2, CO, Li3} 
for other proofs and generalizations).  In other words, 
$\RR_1(\A)$ is the union of a subspace arrangement in $\C^n$.  
As shown by Libgober and Yuzvinsky in \cite{LY},   
each subspace of  $\RR_1(\A)$ has dimension at least $2$,   
two distinct subspaces meet only at $\b{0}$, and  
$\RR_d(\A)$ is the union of those subspaces of dimension 
greater than $d$. 

A purely combinatorial description of $\RR_1(\A)$ was given 
in \cite{Fa, LY}. 
A partition $\mathsf{P}=(\mathsf{p}_1\, | \cdots |\,\mathsf{p}_q)$
of $\LL_1(\A)$ is called {\em neighborly} if 
$\left( \abs{\mathsf{p}_j\cap I}\ge \abs{I}-1 \right)$ implies 
$I\subset \mathsf{p}_j$, for all $I\in \LL_2(\A)$. 
To a neighborly partition $\mathsf{P}$, there corresponds 
an irreducible subvariety of $\RR_1(\A)$, 
\begin{equation} 
\label{eq:lp}
L_{\PP}=\Delta_n \cap 
\bigcap_{\{I\in \LL_2(\A) \mid I\not\subset \mathsf{p}_j, \forall j\} }\,
\{\lambda \mid \sum_{i\in I}\lambda_i = 0\}.
\end{equation}
Moreover, $\dim L_{\PP}>0$ if and only if a certain bilinear 
form associated to $\mathsf{P}$ is degenerate. 
Conversely, all components of $\RR_1(\A)$ arise from 
neighborly partitions of sub-arrange\-ments of $\A$.  

In particular, for each $I\in \LL_2(\A)$ with $\abs{I}\ge 3$, 
there is a {\em local} component, 
$L_{I}=\Delta_n\cap \{\lambda \mid \lambda_i=0 \text{ for } i\notin I\}$, 
corresponding to the partition $(I)$ of $\A_I=\{H_i\mid i\in I\}$.
Note that $\dim L_I=\abs{I}-1$, and thus $L_I\subset \RR_{\abs{I}-2}(\A)$.  
At the other extreme, a component of $\RR_1(\A)$ that does not 
correspond to any proper sub-arrangement $\A'$ is called {\em essential}. 
The smallest arrangement for which the resonance variety has an  
essential component (discovered by Falk~\cite{Fa}) 
is the braid arrangement, $\B$, 
with defining polynomial $Q_{\B}=xyz(x-y)(x-z)(y-z)$, 
see Example~\ref{ex:braid}.   

The resonance varieties $\RR_{d}(G,\K)$, 
where $G=G(\A)$ is an arrangement group, and $\K$ is a field 
of positive characteristic, are much less understood.  
The study of the varieties $\RR_{d}(G,\F_{q})$ was started 
in \cite{MS2}.  Recall that the resonance variety $\RR_{d}(\A)$ 
has integral equations, so we may consider its reduction mod~$q$.
As it turns out, there are arrangements $\A$ for which 
$\RR_{d}(G,\F_q)$ does not coincide with $\RR_{d}(\A)$ mod~$q$,
at certain ``exceptional" primes $q$.     
The two varieties have the same number of local components. 
On the other hand, $\RR_d(G,\F_q)$ may have non-local components, 
even though $\RR_d(\A)$ mod~$q$ has none (see \cite{MS2}  and 
Examples \ref{ex:diamond} and \ref{ex:maclane} below).

\begin{question}
Let $G=G(\A)$ be an arrangement group, and $\K$ a field of positive characteristic.
Are all the irreducible components of $R_d(G,\K)$ linear?
\end{question}

\subsection{Characteristic varieties}
\label{subsec:cv}

Let $G$ be a finitely-presented group. 
 For simplicity, assume again that $G$ has torsion-free 
abelianization.  Set $n=b_1(G)$, and fix 
a basis $t_1,\dots,t_n$ for $H_1(G)\cong\Z^n$. 
Let $\K$ be a field, with multiplicative group of units $\K^*$,   
and let $\Hom(G,\K^*)=\Hom(H_1(G),\K^*)\cong{\K^*}^n$ be 
the group of $\K$-valued characters of $G$.  
The {\em $d$-th characteristic variety} of $G$, 
with coefficients in $\K$, is 
\begin{equation} 
\label{eq:vd}
V_d(G,\K)=\{\bo{t}\in {\K^{*}}^n \mid 
\dim_{\K} H^1(G,\K_{\mathbf{t}}) \ge d\},
\end{equation}
where $\K_{\b{t}}$ is the $G$-module $\K$ 
given by the representation $G\xrightarrow{\ab}\Z^n\xrightarrow{\b{t}}\K^*$. 

The characteristic varieties of $G$ form a descending filtration,  
${\K^{*}}^n=V_0\supseteq V_{1}\supseteq \cdots \supseteq V_{n-1}\supseteq V_{n}$, 
which depends only on the isomorphism type of $G$, up to a monomial 
change of basis in the algebraic torus ${\K^*}^n$ (see e.g.~\cite{MS1}).   

Define the {\em depth} of a character $\mathbf{t}:G\to \K^*$ (relative
to the stratification of $\Hom(G,\K^*)={\K^*}^n$ by the characteristic
varieties) to be:
\begin{equation} 
\label{eq:depth}
\dep_{\K}(\b{t})=\max\, \{ d \mid \b{t}\in V_d(G,\K) \}.
\end{equation}
With this notation, 
$V_d(G,\K) \setminus V_{d+1}(G,\K)=\{\b{t}\in {\K^*}^n\mid \dep_{K}(\b{t}) =d\}$.  

The characteristic varieties of $G$ may be identified 
with the determinantal varieties of the Alexander matrix. 
Indeed, let $A_G(\b{t})$ be the matrix $A_G$ evaluated at $\b{t}\in {\K^{*}}^n$.   
For $0\le d< n$, we have:
\begin{equation} 
\label{eq:cval}
V_d(G,\K)=\{\bo{t}\in {\K^*}^n \mid \rank_{\K} A_G(\bo{t}) < n-d\}. 
\end{equation}
Proofs of this result can be found in~\cite{Hi2, Li2} for $\K=\C$, and \cite{MS3}  
for arbitrary $\K$.  By a well-known result from commutative algebra 
(see \cite[pp.~511--513]{Ei}) the $d$-th characteristic variety of $G$ 
is defined by the ideal $\ann\big(\bigwedge^d(\coker A_G)\big)$, the annihilator of the 
$d$-th exterior power of the Alexander module of $G$.

\subsection{Characteristic varieties of arrangements}
\label{subsec:cvarr}

For hyperplane arrangement groups, the characteristic varieties 
are fairly well understood (at least if $\K=\C$).  This is due, in great 
part, to the foundational work by Green and Lazarsfeld~\cite{GL}, Simpson~\cite{Si}, and
Arapura~\cite{Ar} on the structure of the cohomology support loci for local 
systems on quasi-projective algebraic varieties.   We summarize 
the present state of knowledge, as follows. 

\begin{theorem}  
\label{thm:cvarr}
Let $\A$ be an arrangement of $n$ hyperplanes, with group $G=\pi_1(X(\A))$. 
Let $V_{d}(\A)=V_{d}(G,\C)$, $1\le d\le n$, be the characteristic varieties 
of the arrangement.  Then:  

\begin{alphenum}
\item
\label{a}
The components of $V_{d}(\A)$ are subtori of the character torus 
${\C^*}^{n}$, possibly translated by roots of unity, \cite{Ar}.  
Translated subtori indeed do occur, \cite{Su}. 

\item
\label{b}  
The tangent cone at $\mathbf{1}$ to $V_{d}(\A)$ coincides 
with $R_d(\A)$, \cite{CScv, Li2, CO, Li3}.
\end{alphenum}

\end{theorem}

The components of $V_d(\A)$ passing through $\b{1}$ are therefore 
combinatorially determined.  Given a neighborly partition $\PP$ 
of a sub-arrangement $\A'\subset \A$, we will denote by $C_P$ the 
corresponding subtorus of ${\C^*}^n$, so that 
$\operatorname{T}_{\b{1}} C_{\PP}=L_{\PP}$. In particular, 
if $I\in \LL_2(\A)$, we have a local component, $C_I\subset V_{\abs{I}-2}(\A)$.  

It is not known whether the other components of $V_d(\A)$ 
are combinatorially determined.
The non-Fano arrangement has an isolated point in 
$V_2$, see \cite{CScv} and Example~\ref{ex:diamond}. 
The deleted $\operatorname{B}_3$ arrangement 
has a $1$-dimensional component in $V_1$ 
which does not pass through $\b{1}$, see \cite{Su} 
and Example~\ref{ex:deletedB3}. 

We do not know whether Theorem~\ref{thm:cvarr}\eqref{a} holds for $V_d(G,\K)$, 
if $\ch\K>0$.  Clearly,  Theorem~\ref{thm:cvarr}\eqref{b} fails over $\K=\F_2$, 
since $V_1(G,\F_2)={\F_2^*}^n=\{\b{1}\}$, whereas 
$\RR_1(G,\F_2)$ contains local components of the form $L_I$ mod~$2$. 
It also fails (in a more subtle way) 
over $\K=\F_3$, as in Example \ref{ex:maclane}, where $V_1(G,\F_3)$ has 
$8$ components (all local), whereas $R_1(G,\F_3)$ has $9$ components 
(one non-local).

\begin{example}
\label{ex:cvtoy}
For the near-pencil in Examples~\ref{ex:toy}, 
\ref{ex:alextoy}, we have:
\begin{align*}
\RR_1(G,\K)&= \{(\l_1,\l_2,\l_3,\l_4\}\in \K^4 \mid 
\l_1+\l_2+\l_4=0,\  \l_3=0\},\\
V_1(G,\K)&= \{(t_1,t_2,t_3,t_4\}\in {\K^*}^4 \mid t_1t_2t_4=1,\ t_3=1\}.
\end{align*}
Notice that $R_1=L_{124}$, which is the local component 
corresponding to the triple point $v_3=l_1\cap l_2 \cap l_4$. 
Similarly, $V_1=C_{124}$. 
\end{example}


\section{Homology of finite-index subgroups and torsion points on varieties}
\label{sec:fincov}

\subsection{Mod~$q$ first Betti numbers}
\label{subsec:betti}
  
Let $G$ be a finitely-presented group.  For simplicity, 
we will again assume that $H_1(G)$ is torsion-free, 
say of rank $n$.   For a prime $q$, let  
$b_1^{(q)}(G)=\dim_{\,\F_q}H_1(G; \F_q)$ 
be the ``mod $q$ first Betti~number" of $G$.  
Also let $b_1^{(0)}(G)=b_1(G)$ be the usual Betti number.   
Since homology commutes with direct sums, we have
$b_1^{(q)}(G)=\dim_{\,\K}H_1(G; \K)$, for any field
$\K$ of characteristic $q$.  

Let $K\le G$ be a subgroup of finite index $k=\abs{G:K}$. 
Classically, one computes the homology of $K$ by  
Shapiro's lemma:  $H_*(K,\Z)=H_*(G,\Z[G/K])$.  
In \cite{Fx1}, Fox showed how to compute $H_1(K,\Z)$ 
from the invariant factors of $J_G^{\sigma}$, the Jacobian 
matrix of $G$, followed by the permutation representation 
$\sigma:G\to \GL(k,\Z)$ on the cosets $G/K$. 
The only disadvantage of Fox's method is that the matrix  
$J_G^{\sigma}$ may be too large for practical computations.  
The following result from \cite{MS3} refines Fox's method, 
computing the mod~$q$ first Betti numbers of $K$ from the ranks of 
much smaller matrices, provided $K$ is a {\em normal} subgroup, 
and $q\nmid k$.  

\begin{theorem}[\cite{MS3}]  
\label{thm:genlibsak}
Let $G$ be a finitely-presented group, with $H_1(G)\cong \Z^n$. 
Let $\gamma:G\surj\Gamma$ be a homomorphism onto a finite group 
$\Gamma$, with kernel $K_{\gamma}$. Finally, let $\K$ be a field 
such that $q=\ch\K$ does not divide  
the order of $\G$, and such that $\K$ contains all the 
roots of unity of order equal to the exponent of $\G$.  
Then:
\begin{equation*}
\label{eq:genlibsak}
b_1^{(q)}(K_{\gamma}) = n +\sum_{\rho\ne \triv}
n_\rho( \corank J_G^{\rho\circ \gamma}- n_\rho),
\end{equation*}
where the sum is over all non-trivial, irreducible representations 
of $\G$ over the field $\K$, and where $n_{\rho}$ 
is the degree of such a representation $\rho:\G\to \GL(n_{\rho},\K)$.    
\end{theorem}

Since $\C$ is algebraically closed, and every complex, irreducible 
representation of a finite abelian group has degree $1$, 
the next corollary follows at once.   

\begin{cor}[Libgober~\cite{Li2}, Sakuma~\cite{Sa}, Hironaka~\cite{Hi2}]
\label{cor:libsak}
If $\gamma:G\surj\Gamma$ is a homomorphism of $G$ onto a 
finite abelian group $\G$, then:
\begin{equation} 
\label{eq:b1k}
b_1(K_{\gamma}) = n +\sum_{\triv\ne\rho\in \Hom(\G,\C^*)}
( \corank J_G^{\rho\circ \gamma}- 1).
\end{equation}
\end{cor}

\subsection{Betti numbers and torsion points}
\label{subsec:bettitors}
If $\Gamma$ is a group of prime order $p$, we may reinterpret  
the formula from Theorem~\ref{thm:genlibsak}, in terms of  
depth of $p$-torsion points on the character variety of $G$ 
(relative to the stratification by the characteristic varieties).     
More precisely, let $q=0$, or $q$ a prime, $q\ne p$.  
Take a field $\K$ of characteristic $q$  
which contains all $p$-roots of unity---for example,
$\K=\C$, or $\K=\F_{q^s}$, where $s=\ord_p(q)$ 
is the least positive integer such that $p \mid (q^s-1)$.  
Viewing $\Gamma=\Z_p$ as the subgroup of $p$-roots of 
unity in $\K^*$, we may identify $\Hom(G,\Z_p)$ with 
the subgroup of $p$-torsion points on $\Hom(G,\K^*)$. 

\begin{cor}[\cite{MS3}]
\label{cor:torsdepth}
If $K=\ker (\gamma: G\surj \Z_p)$, then: 
\begin{equation*} 
\label{eq:tordep}
b_1^{(q)}(K) = n +(p-1)\dep_{\K}(\gamma).
\end{equation*}
\end{cor}

In view of this Corollary, we are led to define the following 
invariants of the group $G$.  Let 
\begin{equation}
\label{eq:betainv}
\beta_{p,d}^{(q)}(G) =  \#\big\lbrace K \lhd G \mid [G:K]=p 
\,\text{ and }\, b_1^{(q)}(K) = b_1^{(q)}(G) + (p-1) d \,\big\rbrace 
\end{equation}
be the number of index~$p$,
normal subgroups of $G$ for which the mod~$q$ first Betti number
jumps by $(p-1)d$, when compared to that of $G$. 
Note that $\sum_{d\ge 0} \beta_{p,d}^{(q)}(G)=\frac{p^n-1}{p-1}$.

Now let 
\begin{equation}
\label{eq:torsn}
\Tors_N({\K^*}^n)=\{\b{t}\in {\K^*}^n\mid 
\mathbf{t}^N=\mathbf{1}\ \text{and}\ \mathbf{t}\ne \mathbf{1}\}
\end{equation}
be the set of points on the algebraic torus ${\K^*}^n$ of order 
precisely equal to $N$.  For a subset $V\subset {\K^*}^n$, put  
$\Tors_{N}(V)=\Tors_N({\K^*}^n)\cap V$.  
With these notations, the $\beta$-invariants of $G$ may 
be computed from the characteristic varieties, as follows.  

\begin{theorem}[\cite{MS3}]
\label{thm:torscount}
Let $\K$ be a field of characteristic $q\ne p$, containing
all the $p$-roots of unity.  Then:
\begin{equation*} 
\label{eq:torcnt}
\beta_{p,d}^{(q)}(G)=\frac{\#\Tors_p(V_{d}(G,\K)
\setminus V_{d+1}(G,\K))}{p-1}.
\end{equation*}
\end{theorem}

An analogous formula holds for $q=p$, with the characteristic varieties 
replaced by the resonance varieties over the field $\F_p$. 
Let $G=\langle x_1,\dots , x_n \mid r_1,\dots ,r_m\rangle$ 
be a finitely-presented group, with $H_1(G)\cong\Z^n$ and $H_2(G)$ torsion-free.   
Let $G/\gamma_3 G$ be the second nilpotent quotient of $G$ (see \ref{subsec:LCS}).  
The following invariants of $G$ (in fact, of $G/\gamma_3 G$) 
were introduced in \cite{MS2}:
\begin{equation} 
\label{eq:nupd}
\nu_{p,d}(G)=\#\big\lbrace K \lhd G/\gamma_3 G \mid \abs{G/\gamma_3 G:K}=p
\ \text{and}\ b_1^{(p)}(K) = n + d \:\big\rbrace.
\end{equation}
We then have:
\begin{theorem}[\cite{MS2}]
$
\DS{\nu_{p,d}(G)=\frac{\abs{R_d(G,\F_p)\setminus R_{d+1}(G,\F_p)}}{p-1}}.
$
\end{theorem}


\section{Congruence covers and polynomial periodicity}
\label{sec:cong}

\subsection{Congruence covers}
\label{subsec:congcovers}
Let $X$ be  a finite cell complex (up to homotopy equivalence).  
Assume $H_1(X)\cong \Z^n$.  Attached to $X$ there is a canonical 
sequence of {\em congruence covers}, $\{X_N\}_{N\in \N}$, defined as 
follows:   $X_N$ is the regular $(\Z_N)^n$-cover of $X$ 
determined by the projection 
\begin{equation}
\label{eq:congproj}
\pi_1(X)\xrightarrow{\ab} H_1(X,\Z) \xrightarrow{\bmod N}H_1(X;\Z_N).
\end{equation}
The sequence of congruence covers depends, up to homotopy, 
only on $X$. Indeed, a homotopy equivalence $X\xrightarrow{\simeq} Y$ lifts to 
homotopy equivalences $X_N\xrightarrow{\simeq} Y_N$, for all $N\ge 1$. 
Thus, the sequence of Betti numbers $\{b_1(X_N)\}_{N\in \N}$ is 
a homotopy-type invariant of $X$.  

The first Betti numbers of the congruence covers of $X$   
can be expressed in terms of the characteristic varieties  
of $G=\pi_1(X)$.  Indeed, we may reformulate the results  
of Libgober, Sakuma, and Hironaka (Corollary~\ref{cor:libsak}), 
in this particular case, as follows.
\begin{theorem}
\label{thm:reflibsak}
For every finite cell complex $X$ with $H_1(X)\cong \Z^n$, 
\begin{equation*}
\label{eq:b1cong}
b_1(X_N)=n+\sum_{\b{t}\in \Tors_N({\C^*}^n)} \depth_{\C}(\b{t}).
\end{equation*}
\end{theorem}

As a simple example, take $X=\C\setminus \{\text{$n$ points}\}$.    
Then $\pi_1(X)=F_n$, and  $V_{<n}(F_n,\C)={\C^{*}}^n$,  
$V_n(F_n,\C)=\{ \mathbf{1} \}$.  Thus,  $b_1(X_N)=n+(N^n-1)(n-1)$,
which of course can also be seen by computing Euler characteristics.  

\subsection{Polynomial periodicity}
\label{subsec:polycong}
As first noted by Peter Sarnak (who coined the term), the sequence 
of first Betti numbers of congruence covers of $X$ is 
polynomially periodic. More precisely:

\begin{theorem}[Sarnak-Adams~\cite{SarA}, Sakuma~\cite{Sa}]
\label{thm:sarnak}
For every finite cell complex $X$, 
the sequence $b_1(X_N)$ is {\em polynomially periodic}, i.e., 
there exists an integer $T\ge 1$, 
and polynomials $P_1(x),\dots ,P_{T}(x)$, such that 
\begin{equation*} 
\label{eq:b1xn}
b_1(X_N)=P_i(N), \quad \text{if $N\equiv i\,\bmod T$}.
\end{equation*}
\end{theorem}

We shall call the minimal such $T$ the {\em period} of the 
sequence of Betti numbers.    
Notice that the polynomials $P_{i}$ and the period $T$ 
are homotopy-type invariants of $X$.  
Theorem~\ref{thm:sarnak} readily follows from Theorem~\ref{thm:reflibsak}, 
by means of a deep result of Laurent~\cite{La} (originally 
conjectured by Serge Lang).

\begin{theorem}[Laurent~\cite{La}]
\label{thm:laurent}
If $V$ is a subvariety of ${\C^*}^n$,  
then $\Tors_N(V)=\bigcup_{i=1}^{v} \Tors_N(S_i)$, where 
$S_i$ are subtori of ${\C^*}^n$, possibly translated by 
roots of unity.
\end{theorem}

\begin{remark} 
\label{rem:congbeta}
If $N=p$ is prime, we may express the Betti 
number $b_1(X_p)$ in terms of the $\beta$-invariants 
of $G=\pi_1(X)$:
\begin{equation}
\label{eq:b1xp}
b_1(X_p)=n+(p-1)\sum_{d\ge 1} d\, \beta_{p,d}^{(0)}(G).
\end{equation}
Thus, if the $\beta$-invariants are given by different 
polynomials, for different primes $p$, non-trivial periodicity 
in the Betti numbers of congruence covers is likely to ensue. 
\end{remark}

\begin{remark}
\label{rem:otherper}
The above analysis can be carried out in the slightly more general 
context considered by Sakuma~\cite{Sa}.  Let $\pi_{A}:H_1(X)\surj A$ 
be a homomorphism from $H_1(X)\cong \Z^n$ onto a finitely-generated 
abelian group $A$. There is then a sequence of {\em $\pi_{A}$-congruence 
covers}, $X_{\pi_{A},N}$, defined by the projection 
\begin{equation}
\label{eq:piA}
\pi_1(X)\xrightarrow{\ab} H_1(X) \xrightarrow{\pi_{A}} A \surj A\otimes \Z_N.
\end{equation}  
Unlike the case discussed above (where $A=\Z^n$) 
this sequence of covers is not necessarily a homotopy-type invariant of $X$. 
Nevertheless, the sequence of Betti numbers, $b_1(X_{\pi_{A},N})$, 
is still polynomially periodic, but the polynomials and the period 
depend also on $\pi_{A}$.   
\end{remark}

\subsection{Congruence covers of arrangements}
\label{subsec:arrcong}  
Let $\A$ be a hyperplane arrangement, with complement $X=X(\A)$.  
Up to homotopy equivalence, $X$ is a finite cell complex.      
Let $X_N=X_N(\A)$ be its congruence covers.   
By Theorem \ref{thm:sarnak}, the sequence $b_1(X_N(\A))$ is 
polynomially periodic.  As far as we know, there are no examples in the 
literature where non-trivial periodicity occurs (i.e., where $T>1$).  
This is probably due to the following observation, which follows 
immediately from  Theorem~\ref{thm:cvarr}\eqref{a} 
and Theorem~\ref{thm:reflibsak}. 

\begin{prop}
\label{prop:period1}
Suppose the intersection of any two components of the characteristic 
varieties $V_d(\A)$ ($d\ge 1$) is the identity.  Then $b_1(X_N(\A))$ 
is a polynomial in $N$ (of degree equal to the maximum dimension of 
any of those components).   
\end{prop}

If $n\le 6$, the above condition is satisfied, and so 
$b_1(X_N(\A))$ is a polynomial sequence. 
If $n\ge 7$, the condition may fail:  two components may meet 
outside the origin, or miss each other altogether.  
This leads to the following result. 

\begin{prop}
\label{prop:arrper}
There exist hyperplane arrangements $\A$ for which the sequence 
$b_1(X_N(\A))$ is polynomially periodic, with period $T>1$. 
\end{prop}

In Section~\ref{sec:examples}, we give several such examples.   
The simplest one is the non-Fano plane, with $n=7$ (Example~\ref{ex:diamond}). 
The Betti numbers $b_1(X_N(\A))$ are given by two polynomials, 
$P_1(N)=9N^2-3$ and $P_2(N)=9N^2-2$, depending on whether  
$N$ is even or odd  (thus, the period is $T=2$). 
In all the examples, the non-trivial periodicity of the Betti numbers 
can be traced to the presence 
of isolated points of order $2$ in $V_2(\A)$.  It can also 
be explained by the dependence of $\beta_{p,d}^{(0)}$ (viewed as 
a polynomial in $p$) on the parity of the prime $p$ 
(see Remark~\ref{rem:congbeta}). 

\begin{remark}
\label{rem:co}
Similar results may be obtained for the $\pi_{A}$-congruence covers 
of $X=X(\A)$ mentioned in Remark~\ref{rem:otherper}.  One such family 
of covers was recently considered by Cohen and Orlik~\cite{CO2}.  
They take $A=\Z$, and $\pi_A=(1,\dots,1)^{\top}:\Z^n \to \Z$, where $n=\abs{\A}$, 
and study the resulting sequence of $N$-fold cyclic covers, 
$X_{(N)}=X_{\pi_{A},N}$ (the motivation being that, if $X^*=X(\b{d}\A)$, 
then $X^*_{(n)}$ is the fiber of the Milnor fibration, $Q_{\A}: X(\A)\to \C^*$).  
Again, non-trivial polynomial periodicity appears in the sequence 
$b_1(X_{(N)})$.  But the non-trivial periodicity of the sequence of Betti numbers 
of congruence covers is a much rarer phenomenon than that of cyclic covers. 
For example, if $X=X(\B)$ is the complement of the braid arrangement, 
the sequence $b_1(X_{(N)})$ has period $T=3$, whereas $b_1(X_N)$ has 
trivial period ($T=1$).  
\end{remark}


\section{Hirzebruch covering surfaces}
\label{sec:hirzebruch}

\subsection{Branched covers}
\label{subsec:def}

Let $\A=\{H_1,\dots,H_n\}$ be a central arrangement of planes 
in $\C^3$, with complement $X(\A)$.  
Let $\ov{\A}=\{\ov{H}_1,\dots,\ov{H}_n\}$ 
be the corresponding arrangement of projective 
lines in $\CP^2$, with complement 
$X(\ov{\A})=\CP^2\setminus \bigcup_{i=1}^n \ov{H}_i$. 
Notice that $X(\ov{\A})$ is diffeomorphic to the 
complement of any decone of $\A$.  In particular, 
$H_1(X(\ov{\A}))=\Z^{n-1}$.  

For each plane arrangement $\A$, and positive integer $N$, 
Hirzebruch defined in~\cite{Hz} a compact, smooth, complex algebraic surface, 
$M_N(\A)$, as follows.  Let $X_N(\ov{\A})$ be the congruence 
cover of $X(\ov{\A})$, and let $\widehat{X}_N(\ov{\A})$ be 
the associated branched cover (or, ``Kummer cover") of $\CP^2$.  
This is an algebraic surface with normal singularities 
which ramifies over the projective plane, with the 
arrangement as the branching locus.  
By definition, the {\em Hirzebruch covering surface} $M_N(\A)$ 
is the minimal desingularization of $\widehat{X}_N(\ov{\A})$.  
The resolution of singularities can be described concretely:  
If a point $z\in \CP^2$ 
lies on $r$ lines of $\ov{\A}$, then the preimage of $z$ in 
$\widehat{X}_N(\ov{\A})$ consists of $N^{n-1-r}$ points, 
forming an orbit of the Galois group $(\Z_N)^{n-1}$. 
For $r\ge 3$, each such point $\hat{z}$ is singular. 
The minimal resolution of $\hat{z}$ replaces it by 
a smooth curve of genus $2^{r-3}(r-4)+1$ and self-intersection 
$-2^{r-2}$.

\subsection{Chern numbers}
\label{subsec:chern} 

In \cite{Hz}, Hirzebruch computed the Chern numbers 
of $M=M_N(\A)$, in terms of simple 
combinatorial data associated to $\A$. 
The algebraic surface $M$ has Chern numbers  
$c_1^2$ (equal to the self-intersection number 
of a canonical divisor) and $c_2$ 
(equal to the Euler characteristic of $M$).  
If $\A_n$ is an arrangement for which $\ov{\A}_n$ 
is a pencil of $n$ projective lines, then $c_1^2$ and $c_2$ 
depend only on $n$ and $N$, see Example~\ref{ex:pencil}. 
Otherwise, let $n=\abs{\A}$, $s=\abs{\LL_2(\A)}$, and 
$m_r=\#\{ I\in \LL_2(\A) \mid \abs{I}=r\}$.    
Recall that $\sum_{r} m_r (r-1)=b_2$, where $b_2=b_2(X(\A))$.    
\begin{theorem}[Hirzebruch~\cite{Hz}]
\label{thm:hirzebruch}
The Chern numbers of $M_N(\A)$ are given by
\begin{align*}
c_1^2&=\big( (3b_2-s-5n+9)N^{2} - 4(b_2-n) N  + (b_2+n+m_2) \big) N^{n-3},\\
c_2&=\big( (b_2-2n+3)N^{2} - 2(b_2-n) N  + (b_2+s-m_2) \big) N^{n-3}.
\end{align*}
\end{theorem}

For example, if $\B$ is the braid arrangement,  
with $n=6$, $s=7$, $m_2=3$, and $b_2=11$,   
then $M_2(\B)$ has $c_1^2=0$ and $c_2=24$ 
(in fact, as noted by Hirzebruch, $M_2(\B)$ is the Kummer surface). 

\subsection{Irregularity}
\label{subsec:bettihirz}
The Betti numbers of the compact algebraic surface $M=M_N(\A)$ 
can be computed from $c_2=\chi(M)$ by Poincar\'{e} duality, 
once $b_1(M)$ is known. 
The first Betti number $b_1(M)$---equal to twice the 
irregularity $q(M)$---depends on more subtle data than 
the Chern numbers.  Several methods for computing 
$b_1(M_N(\A))$ were given by Ishida \cite{Ish}, 
Gl\"{a}ser~\cite{Gl}, Zuo~\cite{Zu}, and Hironaka~\cite{Hi1}.  
For example, Ishida found $b_1(M_5(\B))=60$; more generally, 
Gl\"{a}ser found $b_1(M_N(\B))=5(N-1)(N-2)$.  
Hironaka produced detailed tables of computations of 
$b_1$ (and also $c_1^2$, $c_2$) for surfaces $M_N(\A)$ with 
$n\le 7$ and $N\le 5$.  

In \cite{Sa}, Sakuma gave a new and powerful formula  
for computing the first Betti numbers of finite abelian 
branched covers, in a wider context, in terms of nullities 
of Alexander matrices.  
We find it convenient to reinterpret Sakuma's formula in terms 
of depth of complex characters. We do this only for 
Hirzebruch covering surfaces, though a similar interpretation     
in a more general context is possible.  

Let $x_1,\dots ,x_n$ be the meridians of $\A$, which recall 
generate $H_1(G(\A))\cong \Z^n$.  
For a character $\b{t}\in\Hom(H_1(G(\A)),\C^*)\cong {\C^*}^n$, 
and a sub-arrangement $\A'\subset \A$, let $\b{t}|_{\A'}$ 
be the pull-back of $\b{t}$ to $\Hom(H_1(G(\A')),\C^*)$, 
via the obvious inclusion $H_1(G(\A'))\to H_1(G(\A))$.
Also, let 
$\A_{\b{t}}=\{H_i\in \A \mid \b{t}(x_i)\ne 1\}$ be the sub-arrange\-ment 
on which $\b{t}$ is supported. 
Finally, recall that $\dep_{\C}(\b{t})=\max\, \{ d \mid \b{t}\in V_d(\A) \}$.  
With this notation, Sakuma's formula \cite{Sa} can be written as follows.

\begin{theorem}
\label{thm:sakuma}
The first Betti number of $M_N(\A)$ is given by
\[
b_1(M_N(\A))=\sum_{\b{t}\in \Tors_N({\C^*}^n)} 
\depth_{\C}(\b{t}|_{\A_{\b{t}}}).
\]
\end{theorem}

For example, the characteristic variety of the braid arrangement, 
$V_1(\B)$, has five $2$-dimensional components, one of which is non-local, 
see Example~\ref{ex:braid}.  Applying Theorem~\ref{thm:sakuma}, 
one readily recovers Gl\"{a}ser's formula.\footnote{%
In hindsight, Falk's discovery of the fifth component 
of $\RR_1(\B)$---the tangent cone at the origin to $V_1(\B)$---was 
predicted 14 years earlier by Ishida and Gl\"{a}ser!} 

In \cite{Ta}, Tayama used Sakuma's formula (in its original form) 
to complete Hironaka's tables of Betti numbers of $M_N(\A)$, also for 
$n\le 7$, but for arbitrary $N$.  
Moreover, he computed the Betti numbers for pencils of 
arbitrary size, finding:
\begin{equation}
\label{eq:proj}
b_1(M_N(\A_n)):=b(N,n)=(N-1)\Big((n-2)N^{n-2}-2\sum_{k=0}^{n-3} N^k\Big).
\end{equation}
As a consequence, Tayama derived the lower bound
\begin{equation}
\label{eq:tayama}
b_1(M_N(\A))\ge \sum_{r\ge 3} m_r b(N,r) + \beta(\A) b(N,3),
\end{equation}
where $\beta(\A)$ is the number of braid sub-arrangements of $\A$.  

\begin{remark}
\label{rem:tayama}
An inspection of the tables in \cite{Ta} shows that 
equality holds in \eqref{eq:tayama} for all arrangements $\A$ 
with $n\le 7$.  This is not surprising, since, in that range, the only 
sub-arrangement of $\A$ that can have an essential component 
in its characteristic varieties is the braid arrangement $\B$.  

If $n\ge 8$, though, Tayama's inequality can very well be strict, 
as the arrangements in Examples \ref{ex:deletedB3}, \ref{ex:B3}, 
\ref{ex:pappus}, and \ref{ex:ziegler} show.   For each such 
arrangement $\A$, the characteristic varieties have essential 
components that do not come from braid sub-arrangements, 
and some of the torsion points on these components add 
to the Betti number $b_1(M_N(\A))$. 
\end{remark}

\subsection{Polynomial periodicity}
\label{subsec:polyhirz} 

By Theorem~\ref{thm:sarnak}, the sequence of first Betti numbers 
of (unramified) congruence covers, $b_1(X_N(\A))$, is polynomially 
periodic.  It it by no means evident, but still true, that the same 
holds for the ramified congruence covers, and their minimal 
desingularizations. 
\begin{theorem}[Hironaka \cite{Hi0}, Sakuma~\cite{Sa}]
\label{thm:hirsak}
For every arrangement $\A$, the sequence $b_1(M_N(\A))$ 
is polynomially periodic.
\end{theorem}

If $n\le 7$, the sequence $b_1(M_N(\A))$ is in fact given 
by a polynomial in $N$, see Remark~\ref{rem:tayama}.  According to  
Hironaka~\cite{Hi3}, the first (and seemingly the only) example of 
an arrangement for which non-trivial periodicity occurs was given by Zuo 
in \cite{Zu}.  This is the Hessian arrangement, $\mathcal{H}$, 
with defining polynomial 
$Q_{\mathcal{H}}=xyz\prod_{i=0}^{2} \prod_{j=0}^{2} (x+\omega^i y+\omega^j z)$,  
where $\omega=e^{2\pi\ii/3}$. 
The Betti numbers $b_1(M_N(\mathcal{H}))$ are given by two distinct  
polynomials in $N$, depending on 
whether $N$ satisfies certain divisibility 
conditions.  Thus, the period $T$ is at least $2$, 
though the precise value of $T$ does not seem to be known.  

The next result shows that periodicity in $b_1(M_N(\A))$ occurs for the minimum 
possible value of $n=\abs{\A}$.  Moreover, the period can be identified 
precisely. 

\begin{prop}
\label{prop:hirzper}
There exists a hyperplane arrangement $\A$ with $\abs{\A}=8$ 
for which the sequence 
$b_1(M_N(\A))$ is polynomially periodic, with period $T=4$. 
\end{prop}

The arrangement $\A$ is the deleted $\operatorname{B}_3$-arrangement   
(unlike Zuo's example, this is a complexified real arrangement).  
Explicit polynomials describing the periodic behavior 
of $b_1(M_N(\A))$ are given in Example~\ref{ex:deletedB3}

Another instance of non-trivial 
periodicity appears in Example~\ref{ex:ziegler}.  The two Ziegler 
arrangements, $\A_1$ and $\A_2$, have the same multiplicities, 
and so the complex surfaces $M_N(\A_1)$ and $M_N(\A_2)$ have the same Chern numbers.  
Moreover, both sequences $b_1(M_N(\A_i))$ have period $4$, and both start 
with $0$, $108$, and $2,110$, but $b_1(M_4(\A_1))=13,932$ 
and  $b_1(M_4(\A_2))=13,930$.


\section{Counting finite-index subgroups}
\label{sec:counting}

\subsection{Hall formulas}
\label{subsec:hall}
Let $G$ be a finitely-generated group. 
For each positive integer $k$, let
$a_{k}(G)$ be the number of index $k$ subgroups of $G$. 
As noted by Marshall~Hall~\cite{Ha}, these numbers can be 
computed recursively, starting from $a_1(G)=1$, 
by means of the formula
\begin{equation}
\label{eq:hall}
a_{k}(G)=\frac{1}{(k-1)!}\, \abs{\Hom(G,S_k)} -
\sum_{l=1}^{k-1} \frac{1}{(k-l)!}\, \abs{\Hom(G,S_{k-l})}\, a_{l}(G), 
\end{equation}
provided one knows how to count representations from $G$ to 
the symmetric groups.   (See the survey by Lubotzky~\cite{Lu} 
for an overview of the subject.)    

To do this count, it is convenient to introduce  
the Hall invariants of $G$, so named after Philip Hall, who 
first considered them in \cite{HaP}. 
Let $\Gamma$ be a finite group.  The {\em $\Gamma$-Hall invariant} of $G$ 
is the number of surjective representations of $G$ to $\Gamma$,
up to automorphisms of $\Gamma$:
\begin{equation}
\label{eq:delta}
\delta_{\Gamma}(G)=\abs{\Epi(G,\Gamma)/\Aut\Gamma}.
\end{equation}
As noted by P.~Hall, 
\begin{equation}
\label{eq:hom}
\abs{\Hom(G,\Gamma)}=\sum_{H\le\G}\abs{\Aut(H)}\, \delta_{\Gamma}(G).
\end{equation}
By \eqref{eq:hall} and \eqref{eq:hom}, the computation of 
$a_k(G)$ reduces to the computation 
of $\delta_{\Gamma}(G)$, for all subgroups $\Gamma\le S_k$. 
For example, if $H_1(G)\cong \Z^n$, we have:
\begin{align}
\label{a2}
a_2&= 2^{n}-1, 
\\
\label{a3}
a_3&=\tfrac{1}{2}(3^{n}-1) + 3\delta_{S_3},
\\
\label{a4}
a_4&=\tfrac{1}{3}(2^{n+1}-1)(2^n-1) + 
4(\delta_{D_8} + \delta_{A_4} + \delta_{S_4}).
\end{align}

\subsection{Abelian representations}
\label{subsec:abel}
If $\Gamma$ is finite abelian, an explicit formula 
for $\delta_{\Gamma}(G)$ was derived in \cite{MS3}.  
For simplicity, we will give this formula just in the case when 
$G$ has torsion-free abelianization (if $G=F_n$, a similar formula 
was obtained by Kwak, Chun, and Lee in \cite{KCL}). 

Write $\G=\prod_{p \mid \, \abs{\G}} \G_p$, where 
$p$ runs through all the primes dividing the order of $\G$, 
and $\G_p$ is a finite abelian $p$-group. Clearly, 
$\delta_{\G}(G) =\prod_{p \mid \, \abs{\G}} \delta_{\G_p}(G)$,    
and so it is enough to compute $\delta_{\G_p}(G)$.

\begin{theorem}[\cite{MS3}]
\label{thm:delabel}  
Let $G$ be a finitely-generated group, with $H_1(G)\cong \Z^n$.  
If $\G_p=\bigoplus_{i=1}^{k}\Z_{p^{\nu_i}}$ 
is an abelian $p$-group, then:
\begin{equation*}
\label{eq:delpgroup}
\delta_{\G_p}(G)=
\frac{p^{|\nu|(n-1)-2\langle\nu\rangle} \varphi_n(p^{-1})}
{\varphi_{n-k}(p^{-1})\prod_{r\ge 1} \varphi_{m_r(\nu)} (p^{-1})}, 
\end{equation*}
where $|\nu|=\sum_{i=1}^{k} \nu_i$, 
$\langle\nu\rangle=\sum_{i=1}^{k}(i-1)\nu_i$,  
$m_r(\nu)=\#\{ j\mid \nu_j=r\}$, 
and $\varphi_m(t)=\prod_{i=1}^{m} (1-t^i)$.
\end{theorem}

\subsection{Metabelian representations}
\label{subsec:metaabel}
The next level of difficulty is presented by finite metabelian groups.  
We have a formula in that case, too, but only for a special 
class of (split) metabelian groups, which includes the metacyclic 
groups considered by Fox in \cite{Fx70}.  More subtle information 
about $G$ enters in this formula---not just $n=b_1(G)$, but 
also the stratification of $\Hom(G,\K^*)$ by the characteristic 
varieties $V_d(G,\K)$, over certain Galois fields $\K=\F_{q^s}$. 

For two distinct primes $p$ and $q$, consider the semidirect product  
$\Z_q^s \rtimes_{\s} \Z_p$, where $s=\ord_p(q)$ is the order of 
$q\bmod p$ in $\Z_p^*$, and $\sigma$ is an automorphism of $\Z_q^s$  
of order $p$.  The group $\Z_q^s \rtimes_{\s} \Z_p$ is 
independent of the choice of such $\sigma$. 
Viewing $\Z_q^s$ as the additive group of the
field $\F_{q^s}=\F_q(\xi)$, where $\xi$ is a primitive
$p$-th root of unity, we may take $\sigma=\cdot \xi\in\Aut(\F_q(\xi))$.  

\begin{theorem}[\cite{MS3}]
\label{thm:delmpqs}  
Let $G$ be a finitely-presented group. Then:
\begin{equation*}
\label{eq:dmpq}
\delta_{\Z_q^s \rtimes_{\s} \Z_p}(G) =\frac{p-1}{s(q^s-1)} \sum_{d\ge 1}
\beta_{p,d}^{(q)}(G) (q^{sd}-1).
\end{equation*}
\end{theorem}

The theorem applies to the symmetric group 
$S_3=\Z_3\rtimes_{\sigma} \Z_2$,  
where $\sigma=(-1)$, and to the alternating group 
$A_4=\Z_2^2\rtimes_{\sigma} \Z_3$,   
where $\sigma=\bigl(\begin{smallmatrix}0&1\\1&1\end{smallmatrix}\bigr)$.   
Combined with Theorem~\ref{thm:torscount}, it yields:
\begin{align}
\label{eq:dels3}
\delta_{S_3}(G)&=
\tfrac{1}{2} \sum_{d\ge 1}\:
\abs{\Tors_2(V_{d}(G,\F_3)\setminus V_{d+1}(G,\F_3))}\, (3^{d}-1), 
\\  
\label{eq:dela4}
\delta_{A_4}(G)&=\tfrac{1}{3} \sum_{d\ge 1}\:
\abs{\Tors_3(V_{d}(G,\F_4)\setminus V_{d+1}(G,\F_4))}\,  (4^{d}-1).
\end{align}

Returning now to the low-index subgroups of $G$, we see that $a_3(G)$   
can be computed from \eqref{a3} and \eqref{eq:dels3}.  
As for $a_4(G)$, the term $\delta_{A_4}$ from \eqref{a4} 
can be computed from \eqref{eq:dela4}, 
but the terms $\delta_{D_8}$ and $\delta_{S_4}$ cannot 
be handled by this method. 

\subsection{Counting normal subgroups}
\label{subsec:normal}
Finally, let $a_{k}^{\nor}(G)$ be the number of index $k$, normal
subgroups of $G$.  It is readily seen that
\begin{equation}
\label{eq:aknorm}
a_k^{\nor}(G)=\sum_{\abs{\Gamma}=k} \delta_{\Gamma} (G).
\end{equation}
Using Theorem~\ref{thm:delabel}, we obtain (again, under 
the assumption $H_1(G)\cong \Z^n$):
\begin{align*}
a_2^{\nor}&=2^n-1, 
&a_3^{\nor}&=\tfrac{1}{2}(3^n-1),
\\
a_4^{\nor}&=\tfrac{1}{3}(2^{n+1}-1)(2^n-1),
&a_5^{\nor}&=\tfrac{1}{4}(5^n-1),
\\
a_6^{\nor}&=\tfrac{1}{2}(3^{n}-1)(2^{n}-1)+\delta_{S_3}, 
&a_7^{\nor}&=\tfrac{1}{6}(7^n-1),\\
a_8^{\nor}&=\tfrac{1}{21}(2^{n+2}-1)(2^{n+1}-1)(2^{n}-1) +
\delta_{D_8}+\delta_{Q_8}.  
\end{align*}
In particular, $a^{\nor}_k(G)$ can be computed from the characteristic 
varieties of $G$, if $k\le 7$ (more generally, if $k$ has 
at most two prime factors).  This method breaks down at $k=8$, 
for which $\delta_{D_8}$ and $\delta_{Q_8}$ are also needed.


\section{Lower central series quotients and Chen groups}
\label{sec:lcschen}

\subsection{Lower central series}
\label{subsec:LCS}
Let $G$ be a finitely-generated group.  If $U_1$ and $U_2$ are 
two non-empty subsets of $G$, their commutator subgroup 
is $[U_1,U_2]=\{[u_1,u_2]\mid u_1\in U_1, u_2\in U_2\}$, 
where $[u_1,u_2]=u_1u_2u_1^{-1}u_2^{-1}$.  The commutator 
subgroup of $G$ is $G'=[G,G]$. By repeatedly commuting with 
$G$, a descending series of fully-invariant subgroups is obtained:
\begin{equation}
\label{eq:lcs}
G=\gamma_1 G\ge \gamma_2 G\ge \cdots \ge \gamma_k G\ge \gamma_{k+1} G \ge \cdots 
\end{equation}
where $\gamma_{k+1}G=[\gamma _{k}G, G]$.  This is called the 
{\em lower central series} (LCS) of $G$.     
From the Witt-Hall identities, it follows that   
\begin{equation}
\label{eq:bracket}
[\gamma_{k}G,\gamma_{l}G]\le \gamma_{k+l}G.
\end{equation}
(see \cite[p.~293]{MKS}).  Let $\gr_k G=\gamma_{k}G/\gamma_{k+1}G$ 
be $k$-th lower central series quotient of $G$.  As is well-known, 
$\gr_k G$ is a finitely-generated abelian group (see \cite[p.~295]{MKS}). 
Let $\phi_k(G)=\rank (\gr_k G)$ be its rank.   Notice that $\phi_1(G)=b_1(G)$. 

For the free group $F_n$, the LCS quotients are torsion-free, of 
ranks equal to the Witt numbers:
\begin{equation}
\label{eq:witt}
\phi_k(F_n)=w_k(n):=\tfrac{1}{k}\sum_{d\mid k} \mu(d) n^{k/d},
\end{equation}
where $\mu$ is the M\"{o}bius function, see \cite{MKS}.  
By M\"obius inversion, $\sum_{d\mid k} d\, w_d(n) = n^k$.  
In particular, $\phi_1(F_n)=n$ and $\phi_p(F_n)=\frac{n^p-n}{p}$, for $p$ prime.  
Moreover, $\phi_k(F_n)\le n^k$.  

Now assume $G$ has a presentation with $n$ generators---for example, 
$G$ is the group of an arrangement of $n$ hyperplanes.  The 
projection $F_n\surj G$ induces surjections  $\gr_k F_n \surj \gr_k G$.  
Thus, $\phi_k(G)\le \phi_k(F_n) \le n^k$.  A  sharper upper bound was recently 
found by  Newman, Schneider, and Shalev \cite{NSS}:   If $G\ne F_n$, 
there is an $\epsilon>0$ such that 
\begin{equation}
\label{eq:nss}
\phi_k(G) \le (n-\epsilon)^k,\quad  \text{for $k$ sufficiently large}.
\end{equation}

\subsection{Chen Groups}
\label{subsec:Chen}
Let $G/G''$ be the quotient of $G$ by its second derived subgroup.   
The group $G/G''$ is metabelian (in fact, maximal among all metabelian 
quotients of $G$), and, of course, finitely-generated. 
It fits into the exact sequence 
\begin{equation}
\label{eq:ext}
0\to G'/G'' \to G/G'' \to G/G' \to 0.
\end{equation}
The {\em $k$-th Chen group} of $G$ is, by definition, $\gr_k (G/G'')$,  
see \cite{Mu, Ma}.   Let $\theta_k(G)=\phi_k(G/G'')$ be its rank.  

The projection  $G\surj G/G''$ 
induces surjections $\gr_k G\surj \gr_k G/G''$.  Thus, 
\begin{equation}
\label{eq:pt}
\phi_k\ge \theta_k, \quad \text{for all $k$.}
\end{equation}
From \eqref{eq:bracket}, we have  
$G''=[\gamma_2 G,\gamma_2 G]\le \gamma_4 G$, 
and so $\phi_k=\theta_k$ for $k\le 3$. The inequality \eqref{eq:pt} 
is usually strict for $k>3$.  For example (see \cite{Mu, MT, CSai}):
\begin{equation}
\label{eq:chenfn}
\theta_k(F_n)=\binom{n+k-2}{k} (k-1), \quad \text{for $k\ge 2$},
\end{equation}
and this sequence grows much more slowly than the sequence \eqref{eq:witt} 
of Witt numbers.

Now assume $G/G'\cong \Z^n$, with generators $t_1,\dots ,t_n$.  
The Chen groups of $G$ can be determined from the Alexander 
invariant $B=G'/G''$ (viewed as a module over 
$\Z[G/G']\cong \Z[t_1^{\pm 1},\dots, t_n^{\pm 1}]$).  
Indeed, let $\gr_k B=\I^{k}B/\I^{k+1}B$, where $\I$ 
is the augmentation ideal.  Then 
$\gr_k (G/G'')= \gr_{k-2} B$, for $k\ge 2$ (see Massey \cite{Ma}). 
In particular, we have: 
\begin{equation}
\label{eq:hilb}
\sum_{k\ge 0} \theta_{k+2}\, t^k = \Hilb(\gr B),
\end{equation}
where $\gr B=\bigoplus _{k\ge 0} \gr_k B$ (viewed as a module 
over $\gr \Z[G/G']\cong \Z[\l_1,\dots , \l_n]$).  
As noted in \cite{CSpn, CSai}, 
a presentation for $\gr B$ can be obtained 
from a presentation for $B$ via the well-known 
Gr\"obner basis algorithm for finding the tangent 
cone to a variety (cf. \cite{CLO}). 

\subsection{Chen groups of arrangements}
\label{subsec:chenarr}
Let $\A$ be a complex hyperplane arrangement, with group $G=G(\A)$. 
A finite presentation for the Alexander invariant $B=G'/G''$ was    
given in \cite{CSai}.   From this presentation, one can compute 
the Chen groups of the arrangement, $\gr_k (G/G'')$, by the method 
outlined above (see also Westlund~\cite{We} for related work on the 
Chen groups of the ``boundary manifold'' of $\A$).  

It seems very likely that the Chen groups of an arrangement 
are combinatorially determined.  We offer a conjecture which 
makes this combinatorial dependence explicit.  First, we need 
some notation.    Let $\RR_1(\A)=\RR_1(G,\C)$ be the resonance 
variety of $\A$.   Recall that 
$\RR_1(\A)=\bigcup_{i=1}^{v} L_i$, with $L_i$ linear subspaces 
of $\C^n$, where $n=\abs{\A}$. For each $r\ge 1$, let 
\begin{equation}
\label{eq:hr}
h_r=\#\{ L_i \mid \dim L_i=r\}
\end{equation}
be the number of components of $\RR_1(\A)$ of dimension $r$.   
Recall also that $m_r=\#\{ I\in \LL_2(\A)\mid \abs{I}=r\}$. 
We know from the discussion in \S\ref{subsec:rvarr} 
that $h_1=0$, and $h_r\ge m_{r+1}$ for all $r$.  
Moreover, the numbers $h_r$ can be computed directly from 
the intersection lattice $\LL(\A)$, by counting neighborly partitions 
$\PP$ of $\A'\subset\A$, and finding $\dim L_{\PP}$ by linear algebra. 

\begin{conj}[Resonance formula for Chen groups] 
\label{conj:chen}
Let $G=G(\A)$ be an arrangement group. Then:   
\[
\theta_k(G) = \sum_{r\ge 2} h_r  \theta_k(F_r), 
\quad\text{for $k\ge 4$.}
\] 
Moreover, the Chen groups $\gr_k G/G''$ are free abelian 
(of rank $\theta_k$), for all $k\ge 1$. 
\end{conj}

In other words, the conjecture asserts that 
$\theta_k(G) = (k-1) \sum_{r\ge 2} h_r  \binom{k+r-2}{k}$, 
for $k$ sufficiently large. 
If we consider only the local components 
$\{L_I \mid I\in \LL_2(\A)\}$ of $\RR_1(\A)$, 
Conjecture~\ref{conj:chen} implies  
\begin{equation}
\label{eq:ccbound}
\theta_k(G) \ge \theta_k^{\operatorname{cc}}(G):=
(k-1)\sum_{r\ge 3} m_r \binom{k+r-3}{k}. 
\end{equation}
This inequality is actually known to hold, see \cite[Corollary~7.2]{CSai} 
for a direct proof.

Conjecture~\ref{conj:chen} grew out of joint work with 
D.~Cohen~\cite{CSpn, CSai}.  An early version 
of the conjecture appeared in \cite{CSpn}.  It stated\footnote{%
Note the similarity with Tayama's 
lower bound \eqref{eq:tayama} for $b_1(M_N(\A))$.}:  
$\theta_k(G) = \theta_k^{\operatorname{cc}}(G) + \beta(\A) (k-1)$.
Translated into  the present language, 
that earlier conjecture involved only the local 
components of $\RR_1(\A)$, and the non-local components corresponding 
to braid sub-arrangements.   A counterexample was given in \cite{CSai}: 
the Pappus arrangement has an essential component in $\RR_1$, 
and this adds to the ranks of the Chen groups 
(see Example~\ref{ex:pappus}).  On the other hand, in all the 
examples we consider in Section~\ref{sec:examples} (and many others), 
Conjecture~\ref{conj:chen} holds.

\subsection{LCS quotients of arrangements}
\label{subsec:lcsarr}

If $\A$ is a hyperplane arrangement, with group $G=G(\A)$, the LCS ranks 
$\phi_k(G)$ are combinatorially determined.  This was proved 
by Falk in \cite{Fa1}, using Sullivan's $1$-minimal models. 
Even so, a precise combinatorial formula for $\phi_k$ (even for $\phi_3$) 
is not known.  Nevertheless, there is a class of arrangements, singled out 
by Falk and Randell in \cite{FR1}, for which a simple (yet rather 
deep) formula obtains.  

An arrangement $\A=\A_{\ell}$ is called {\em fiber-type} 
if its complement sits atop a tower of {\em linear} fibrations,  
\begin{equation}
\label{eq:ftype}
X=X(\A_{\ell})\xrightarrow{p_{\ell}}X(\A_{\ell-1})\to\cdots\to X(\A_2)
\xrightarrow{p_2} X(\A_1)=\C^*,
\end{equation} 
with fibers $p_i^{-1}(\text{point})=\C\setminus \{\text{$d_i$ points}\}$.  
Thus, $X$ is an Eilenberg-MacLane space $K(G,1)$, with 
fundamental group 
\begin{equation}
\label{eq:semiprod}
G=F_{d_{\ell}}\rtimes_{\rho_{\ell}} \cdots 
\rtimes F_{d_{2}}\rtimes_{\rho_2}  F_{d_{1}}.  
\end{equation}
From the linearity of the fibrations, it follows that each monodromy map 
$\rho_i:\pi_1(X(\A_{i-1}))\to \Aut(F_{d_i})$ factors through the pure 
braid group $P_{d_i}$, and thus acts as the identity on homology.  
A spectral sequence argument shows that the Poincar\'e polynomial 
of the complement factors as $P(X,t)=\prod_{i=1}^{\ell} (1+d_{i} t)$. 
In particular, the {\em exponents} of the fiber-type arrangement, 
$\{d_1=1,d_2,\dots,d_{\ell}\}$, are combinatorially determined.   
 
The {\em LCS formula} of Falk and Randell~\cite{FR1} expresses 
the lower central series ranks, $\phi_k=\phi_k(G)$, 
of a fiber-type arrangement group, $G=G(\A)$,  
in terms of the exponents of $\A$:
\begin{equation}
\label{eq:lcsformula}
\prod _{k\ge 1}(1-t^{k})^{\phi _{k}}=
P(X,-t) = \prod _{i=1}^{\ell}(1-d_{i}t).
\end{equation}
Consequently, 
\begin{equation}
\label{eq:addwitt}
\phi_k(G)=\sum_{i=1}^{\ell} w_k (d_i). 
\end{equation}

As shown by Papadima and Yuzvinsky \cite{PY}, 
the converse holds for arrangements in $\C^3$:  
If $\prod _{k\ge 1}(1-t^{k})^{\phi _{k}}=P(X,-t)$, 
then $\A$ is fiber-type.  

The best-known fiber-type arrangement is 
the braid arrangement $\B_{\ell}$, with exponents 
$\{1,2,\dots,\ell-1\}$.  In that case, the LCS formula was first 
proved by Kohno~\cite{Ko}.  The LCS formula has been an object 
of intense investigation ever since. It has been interpreted  
as a consequence of Koszul duality (Shelton-Yuzvinsky \cite{SY}),   
and has been generalized to hypersolvable arrangements 
(Jambu-Papadima \cite{JP}), and formal rational $K(G,1)$ spaces 
(Papadima-Yuzvinsky \cite{PY}). 

We propose here a formula for the LCS ranks of a certain 
class of arrangements.  The formula looks very similar 
to the LCS formula for fiber-type arrangements, with the exponents 
being replaced by the dimensions of the components of the 
resonance variety.

\begin{conj}[Resonance LCS formula] 
\label{conj:reslcs}
Let $G=G(\A)$ be an arrangement group. If $\phi_4(G)=\theta_4(G)$, then:
\[
\phi_k(G) = \sum_{r\ge 2} h_r \phi_k(F_r), 
\quad\text{for $k\ge 4$.}
\]  
Moreover, the LCS quotients $\gr_k G$ are free abelian 
(of rank $\phi_k$), for all $k\ge 1$. 
\end{conj}

The above formula can also be written as 
$\phi_k(G) = \sum_{r\ge 2} h_r w_k(r)$, or 
\begin{equation}
\label{eq:reslcs}
\prod_{k\ge 1}(1-t^{k})^{\phi_{k}}=
\prod_{i=1}^{v}(1-(\dim L_i)\, t).
\end{equation}
where recall $\RR_1(\A)=\bigcup_{i=1}^{v} L_i$ is the decomposition 
of the resonance variety of $\A$ into (linear) components.  

The equality $\phi_4=\theta_4$  holds for the $\rm{X}_3$, 
$\rm{X}_2$, Pappus, and non-Pappus arrangements 
(Examples~\ref{ex:x3}, \ref{ex:x2}, \ref{ex:pappus}, \ref{ex:nonpappus}).  
The conjecture is verified in those cases (as well as several others), 
at least for $k$ up to $8$.  

There are many non-LCS arrangements for which $\phi_4>\theta_4$, 
for example, the non-Fano and MacLane arrangements 
(Examples~\ref{ex:diamond}, \ref{ex:maclane}). 
In such a situation, the behavior of the sequence $\{\phi_k\}$ is much more 
complicated (see Peeva~\cite{Pe} for a discussion in the non-Fano case).  
Furthermore, the LCS quotients may have torsion, as illustrated in 
Example~\ref{ex:maclane}.  This raises the following question. 

\begin{question}
\label{quest:torslcs}
Let $G=G(\A)$ be an arrangement group.  
Is the torsion in $\gr_k G$ combinatorially determined?
\end{question}


\section{Guide to examples}
\label{sec:glossary}

\renewcommand{\arraystretch}{1.65}

\begin{table}
\caption{Invariants of arrangements}
\label{tab:laundry} 
\begin{minipage}{\columnwidth}
\begin{center}
\renewcommand{\thempfootnote}{\fnsymbol{mpfootnote}}
\setlength{\tabcolsep}{4pt}
\begin{tabular}[t]{@{}p{1.45in}@{\hspace{1pt}}|@{\hspace{5pt}}p{2.95in}}
\hline
Defining polynomial   
& $Q=Q_{\A}$ 
\\
\hline
Multiplicities 
&$n=\abs{\LL_1(\A)}$, $s=\abs{\LL_2(\A)}$, \newline
$m_r=\#\{ I\in \LL_2(\A)\mid \abs{I}=r\}$
\\
\hline
Poincar\'{e} polynomial
&$P(X,t)=1+b_1 t + b_2 t^2 + b_3 t^3$  
\newline\vspace*{0.01in}
\hspace*{8pt} $\big(b_1=n,\: b_2=\sum_{r} m_r (r-1),\: b_3=b_2-n+1\big)$
\\
\hline
Group & $G\cong G^*\times \Z$, where $G^*=G(\b{d}\A)$
\\
\hline
Characteristic varieties
& $V_d(G,\K)$, for ``generic" $\K$, except when noted 
\\
Resonance varieties
& $R_d(G,\K)=\TC_{\b{1}}(V_d(G,\K))$, except when noted 
\\
\hline
$\nu$- and $\beta$-invariants\footnote{Only non-zero values. 
For $d=0$:   \
$\nu_{p,0}=\tfrac{p^n-1}{p-1}-\sum_{d> 0}\nu_{p,d}$, \
$\beta_{p,0}^{(q)}=\tfrac{p^n-1}{p-1}-\sum_{d> 0}\beta_{p,d}^{(q)}$.\\[-2pt]
}
& $\nu_{p,d}=\tfrac{1}{p-1}\abs{\RR_{d}(G,\Z_p)
\setminus \RR_{d+1}(G,\Z_p)}$
\\
&
$\beta_{p,d}^{(q)}=\tfrac{1}{p-1}\abs{\Tors_p(V_{d}(G,\K)\setminus V_{d+1}(G,\K))}$
\\[-2pt]
&
\hspace*{8pt}\vspace*{0.02in}\mbox{$\big( q=0, \K=\C, \text{ or } 
q \text{ prime}, q\ne p, \K=\F_{q^{\ord_{p}(q)}} \big)$}
\\
\hline
 Congruence covers
& First Betti numbers $b_1(X_N)$ and $b_1(M_N)$ 
\\
&
Chern numbers $c_1^2(M_N)$ and $c_2(M_N)$
\\
\hline
 Hall invariants &
$\delta_{S_3}=\tfrac{1}{2} \sum_{d}
\beta_{2,d}^{(3)} (3^{d}-1)$, 
$\delta_{A_4}=\tfrac{1}{3} \sum_{d}
\beta_{3,d}^{(2)} (4^{d}-1)$ 
\\
Low-index subgroups &  
$a_2=2^{n}-1$, 
$a_3^{\nor}=\tfrac{1}{2}(3^{n}-1)$, 
$a_3=a_3^{\nor} + 3\delta_{S_3}$ 
\\
\hline
 LCS ranks\footnote{Values up to $k=8$.    
The general term, for $k\ge 4$, is also given, when 
known from the LCS formula, or, conjecturally, 
from the Resonance LCS formula.
\\[3pt]
}
&    $\phi_1=n$, $\phi_2=\binom{n}{2}-b_2$,  $\phi_k=\rank \gr_k G$
\\
 Chen ranks
&  
$\theta_1=\phi_1$, $\theta_2=\phi_2$, $\theta_3=\phi_3$, 
$\theta_k=\rank \gr_k (G/G'')$  
\\
\hline
\end{tabular}
\end{center}
\end{minipage}
\end{table}

In the next section, we will illustrate the above discussion with 
a collection of examples of hyperplane arrangements, and some of their   
associated invariants.  In this section, we provide a guide 
to those examples, as well as a brief explanation of the computational 
techniques involved. 

In each example, we start with a central arrangement 
in $\C^3$, given by a defining polynomial $Q_{\A}=f_1\cdots f_n$, 
with $f_i=f_i(x,y,z)$ distinct linear forms.  If $\A$ is the 
complexification of a real arrangement, we show a picture of 
the real part of a decone of $\A$. This is an affine arrangement of lines 
in $\R^2$, with the lines ordered as indicated (the $n$-th line being 
the ``line at infinity").  
Otherwise, we merely show the underlying matroid of $\A$, which keeps 
track of the incidence relations.

Table~\ref{tab:laundry} contains the list of objects 
and invariants that we associate with $\A$.   
Some immediate numerical information is extracted from the 
intersection lattice, $\LL(\A)$.  More refined information 
is encoded in the fundamental group of the complement, $G=\pi_1(X(\A))$, 
and its resonance and characteristic varieties. Those 
varieties are depicted schematically, as posets ranked 
by dimension (indicated on vertical scale), and filtered 
by depth (indicated by color scheme: $\text{solid $=$ component  
of $V_1$}$, $\text{open $=$ component of $V_2$}$, etc). 
Furthermore, the components of $V_i$ and $R_i$ 
are distinguished by their shapes.  For example:
\begin{itemize}
\item[$\blacktriangle = $] local component of dimension $2$, 
\item[$\square = $] local component of dimension $3$, 
\item[$\boxplus = $] local component of dimension $4$, 
\item[$\blacklozenge = $] non-local component from braid sub-arrangement,
\item[$\bigstar =$] non-local component from deleted $\operatorname{B}_3$ 
sub-arrangement, etc. 
\end{itemize}

Table~\ref{tab:types} contains the list of invariants 
under consideration, grouped according to the minimal 
input needed to compute them.  Some of these invariants 
can be computed from very simple combinatorial data, others from 
the resonance/characteristic varieties (over various 
fields), while the remaining ones seem to need another kind 
of input, either of a very sophisticated combinatorial nature, 
or altogether non-combinatorial.  

\begin{table}
\caption{Invariants of arrangements, according 
to input needed to compute them}
\label{tab:types} 
\begin{minipage}{\columnwidth}
\begin{center}
\renewcommand{\thempfootnote}{\fnsymbol{mpfootnote}}
\setlength{\tabcolsep}{4pt}
\begin{tabular}[t]{@{}p{1.45in}@{\hspace{1pt}}|@{\hspace{5pt}}p{2.7in}}
\hline
Multiplicities  
&
$P(X,t)$, $c_1^2(M_N)$,  $c_2(M_N)$, $a_2$, $a_3^{\nor}$, 
$a_4^{\nor}$, $a_5^{\nor}$, $a_7^{\nor}$,\newline
$\phi_1=\theta_1$, $\phi_2=\theta_2$, and $\phi_k$ (if $\A$ fiber-type)
\\
\hline
Resonance varieties  
&
$\nu_{p,d}$, $\theta_k$ (conjecturally), \newline  
$\phi_k$ (conjecturally, if $\phi_4=\theta_4$)
\\
\hline
Characteristic varieties 
&
$\beta_{p,d}^{(q)}$, $b_1(X_N)$,  $b_1(M_N)$,  $\delta_{S_3}$, 
$\delta_{A_4}$, $a_3$, $a_6^{\nor}$ 
\\
\hline
Other? 
&
$\delta_{\G}$ (in general), $a_4$, $a_8^{\nor}$,  $\phi_k$ (in general),
\newline
torsion in $H_1(X_N)$, $H_1(M_N)$, $\gr_k(G)$  
\\
\hline
\end{tabular}
\end{center}
\end{minipage}
\end{table}

\renewcommand{\arraystretch}{1.1}

Finally, a few words about how the computations were done: 
\begin{itemize}
\item  The package {\sl Mathematica~4.0}\/ (from Wolfram Research) 
was used to find the braid monodromy presentation of the fundamental 
group, and presentation matrices for various Alexander modules.  
Also, for standard symbolic computations and organization 
of information.
\\[-6pt]

\item  The algebraic geometry and commutative algebra system 
{\sl Macaulay~2}\/ (by Grayson and Stillman~\cite{GS}) was used 
to compute the ideals of minors of various matrices, and their 
primary decompositions.  Also, to compute the tangent cones 
to modules, and their Hilbert series. 
\\[-6pt]

\item  The group theory package {\sl GAP}\/ (from the GAP group~\cite{gap}) 
was used to manipulate group presentations, find finite-index subgroups, 
and compute their abelianizations.  Also, to compute the lower central 
series quotients, via the nilpotent quotient algorithm. 

\end{itemize}


\section{Examples}
\label{sec:examples}


\begin{example}[Central arrangement in $\C^2$]
\label{ex:pencil} 

The most basic arrangement is a pencil of 
$n$ complex lines in $\C^2$.  This is a generic $2$-slice 
of the arrangement $\A_n$ mentioned in \S\ref{subsec:chern}. 
It is a fiber-type arrangement, with exponents $\{1,n-1\}$.

\begin{figure}[ht]
\subfigure{%
\label{fig:pen3}%
\begin{minipage}[t]{0.25\textwidth}
\setlength{\unitlength}{10pt}
\begin{picture}(4.5,2.7)(-3,-1)
\put(0,0){\line(1,1){4}}
\put(-1,2){\line(1,0){6}}
\put(0,4){\line(1,-1){4}}
\put(4.4,-0.5){\makebox(0,0){$1$}}
\put(5.5,1.95){\makebox(0,0){$2$}}
\put(4.4,4.5){\makebox(0,0){$3$}}
\end{picture}
\end{minipage}
}
\setlength{\unitlength}{0.8cm}
\subfigure{%
\label{fig:cvpen3}%
\begin{minipage}[t]{0.15\textwidth}
\begin{picture}(3,2.7)(-0.5,-0.8)
\xygraph{!{0;<6mm,0mm>:<0mm,12mm>::}
[]*D(3){\b{1}}*\cir<3pt>{}
(-[d]*U(3){C_{123}}*-{\blacktriangle}
}
\end{picture}
\end{minipage}
}
\subfigure{%
\label{fig:pen4}%
\begin{minipage}[t]{0.25\textwidth}
\setlength{\unitlength}{10pt}
\begin{picture}(4.5,2.7)(-3,-0.8)
\put(0,0){\line(1,1){4}}
\put(-1,2){\line(1,0){6}}
\put(0,4){\line(1,-1){4}}
\put(2,-0.5){\line(0,1){5}}
\put(4.4,-0.5){\makebox(0,0){$1$}}
\put(5.5,1.95){\makebox(0,0){$2$}}
\put(4.4,4.5){\makebox(0,0){$3$}}
\put(2,5.2){\makebox(0,0){$4$}}
\end{picture}
\end{minipage}
}
\setlength{\unitlength}{0.8cm}
\subfigure{%
\label{fig:cvpen4}%
\begin{minipage}[t]{0.15\textwidth}
\begin{picture}(3,2.7)(-0.5,-0.3)
\xygraph{!{0;<6mm,0mm>:<0mm,16mm>::}
[]*D(3){\b{1}}*\cir<3pt>{}
(-[d]*U(3){C_{1234}}*{\square}
}
\end{picture}
\end{minipage}
}
\caption{\textsf{Pencils of $3$ and $4$ lines}}
\label{fig:pencil}
\end{figure}

\begin{laundry}

\item 
$Q=(y-x)(x-2x)\cdots (y-nx)$.

\item 
$s=m_n=1$.

\item 
$P(X,t)=(1+t)(1+(n-1)t)$.

\item 
$G^*=F_{n-1}$, $G= F_{n-1}\rtimes_{\a} F_1$, where $\a=A_{12\dots n}$. 

\item 
$V_d(G,\K)=\{\b{t}\in {\K^*}^n\mid t_1\cdots t_n=1\}\cong (\K^*)^{n-1}$, for $d<n-1$
\\[3pt]
$V_{n-1}(G,\K)=\{ \mathbf{1} \}$.

\item 
$\beta_{p,n-2}^{(q)}=
\nu_{p,n-2}=\frac{p^{n-1}-1}{p-1}$.

\item 
$b_1(X_N)=(n-2)N^{n-1}+2$.

\item
$b_1(M_N)=(N-1)\big((n-2)N^{n-2}-2\sum_{k=0}^{n-3} N^k\big)$.

\item
$c_1^2(M_N)=N^{n-3}((n- 4)N -n)((n - 2)N -n)$, \\[3pt]
$c_2(M_N)=2(2-n)N^{n-1}+ N^{n-2}+2n$.

\item 
$\delta_{S_3}=\frac{(2^{n-1}-1)(3^{n-2}-1)}{2}$, 
$\delta_{A_4}=\frac{(3^{n-1}-1)(4^{n-2}-1)}{6}$. 

\item 
$a_2=2^{n}-1$, 
$a_3^{\nor}=\tfrac{1}{2}(3^{n}-1)$, 
$a_3=3(3^{n-2}-1)(2^{n-2}+1)+4$. 

\item 
$\phi_{1}=n$, 
$\phi_2=\frac{(n-1)(n-2)}{2}$, 
$\phi_3=\frac{n(n-1)(n-2)}{3}$, 
$\phi_4=\frac{n(n-1)^2(n-2)}{4}$,  \\[1pt]
$\phi_{k}=w_k(n-1)$. 

\item 
$\theta_{1}=n$, 
$\theta_2=\frac{(n-1)(n-2)}{2}$, 
$\theta_3=\frac{n(n-1)(n-2)}{3}$, 
$\theta_4=\frac{(n+1)n(n-1)(n-2)}{8}$,  \\[1pt]
$\theta_{k}=\binom{n+k-3}{k}(k-1)$. 

\end{laundry}
\end{example}


\begin{example}[$\operatorname{X}_3$ arrangement] 
\label{ex:x3} 

This is the smallest arrangement for which the LCS formula does not hold, 
see Falk and Randell \cite{FR86}.

\begin{figure}[ht]
\subfigure{%
\label{fig:dx3}%
\begin{minipage}[t]{0.3\textwidth}
\setlength{\unitlength}{14pt}
\begin{picture}(4.5,2.5)(-2,0)
\multiput(0,1)(0,2){2}{\line(1,0){4}}
\multiput(1,-0.6)(2,0){2}{\line(0,1){4.6}}
\put(0.34,4){\line(2,-3){3.05}}
\put(3.9,-0.55){\makebox(0,0){$1$}}
\put(4.5,1){\makebox(0,0){$2$}}
\put(4.5,3){\makebox(0,0){$3$}}
\put(3,4.5){\makebox(0,0){$4$}}
\put(1,4.5){\makebox(0,0){$5$}}
\end{picture}
\end{minipage}
}
\setlength{\unitlength}{0.8cm}
\subfigure{%
\label{fig:cvx3}%
\begin{minipage}[t]{0.4\textwidth}
\begin{picture}(3,2.5)(-2.5,-0.4)
\xygraph{!{0;<6mm,0mm>:<0mm,16mm>::}
[]*DL(2){\b{1}}*\cir<3pt>{}
(-[dll]*U(3){C_{135}}*-{\blacktriangle}
,-[d]*U(3){C_{236}}*-{\blacktriangle}
,-[drr]*U(3){C_{456}}*-{\blacktriangle}
}
\end{picture}
\end{minipage}
}
\caption{\textsf{The arrangement $\operatorname{X}_3$}}
\label{fig:x3}
\end{figure}

\begin{laundry}

\item 
$Q=xyz(y+z)(x-z)(2x+y)$.

\item 
$n=6$, $s=9$, $m_2=6$, $m_3=3$.

\item 
$P(X,t)=(1+t)(1+5t+7t^2)$.

\item 
$G^*=G(A_{34},A_{24},A_{14},A_{12},A_{135},A_{25}^{A_{35}})$. 

\item 
$V_1(G,\K)=C_{135}\cup C_{236}\cup C_{456}$, 
$V_2(G,\K)=\{ \mathbf{1} \}$. 

\item 
$\beta_{p,1}^{(q)}=\nu_{p,1}=3(p+1)$.

\item 
$b_1(X_N)=3(N^{2}+1)$.

\item 
$b_1(M_N)=3(N-1)(N-2)$.

\item
$c_1^2(M_N)=6N^3(N-2)^2$, 
$c_2(M_N)=3N^3(N^2-4N+5)$.

\item 
$\delta_{S_3}=9$, 
$\delta_{A_4}=12$.  

\item 
$a_2=63$, 
$a_3^{\nor}=364$, 
$a_3=391$. 

\item 
$\phi_{1}=6$, $\phi_2=3$, $\phi_3=6$, $\phi_4=9$, $\phi_5=18$, $\phi_6=27$, 
$\phi_7=90$, $\phi_8=150$.  
Conjecture: $\phi_{k}=3 w_k(2)$. 

\item 
$\theta_{1}=6$, $\theta_2=3$, $\theta_3=6$, $\theta_4=9$,  $\theta_{k}=3(k-1)$. 

\end{laundry}
\end{example}


\begin{example}[Braid arrangement]
\label{ex:braid} 

This arrangement is a generic $3$-slice of $\B_4$.  
It is a fiber-type arrangement, with exponents $\{1,2,3\}$.  
The braid arrangement is the smallest 
arrangement for which the resonance variety contains 
a non-local component (first found by Falk~\cite{Fa}).  

\begin{laundry}

\item 
$Q=xyz(x-y)(x-z)(y-z)$.

\item 
$n=6$, $s=7$, $m_2=3$, $m_3=4$.

\item 
$P(X,t)=(1+t)(1+2t)(1+3t)$.

\item 
$G=P_4$, $G^*=\F_3\rtimes_{\a} \F_2$, where $\a=\{A_{12},A_{13}\}$. 

\item 
$V_1(G,\K)=C_{124}\cup C_{135}\cup C_{236}\cup C_{456} \cup \Pi$, 
where 
\[
\Pi=C_{(16|25|34)}=\{(s,t,(st)^{-1},(st)^{-1},t,s) \mid s,t\in \K^*\}.
\]
$V_{2}(G,\K)=\{\b{1}\}$

\item 
$\beta_{p,1}^{(q)}=\nu_{p,1}=5(p+1)$.

\item 
$b_1(X_N)=5N^{2}+1$.

\item  
$b_1(M_N)=5(N-1)(N-2)$.

\item
$c_1^2(M_N)=5N^3(N-2)^2$, $c_2(M_N)=N^3(2N^2-10N+15)$. 

\item 
$\delta_{S_3}=15$, $\delta_{A_4}=20$. 

\item 
$a_2=63$, 
$a_3^{\nor}=364$, 
$a_3=409$. 

\item 
$\phi_{1}=6$, $\phi_2=4$, $\phi_3=10$, $\phi_4=21$,  $\phi_{k}=w_k(2)+w_k(3)$. 

\item 
$\theta_{1}=6$, $\theta_2=4$, $\theta_3=10$, $\theta_4=15$,  $\theta_{k}=5(k-1)$. 

\end{laundry}

\begin{figure}
\subfigure{%
\label{fig:dbraid}%
\begin{minipage}[t]{0.3\textwidth}
\setlength{\unitlength}{14pt}
\begin{picture}(4.5,2.5)(-1.2,0)
\multiput(0,1)(0,2){2}{\line(1,0){4}}
\multiput(1,0)(2,0){2}{\line(0,1){4}}
\put(0,4){\line(1,-1){4}}
\put(4.5,0){\makebox(0,0){$1$}}
\put(4.5,1){\makebox(0,0){$2$}}
\put(4.5,3){\makebox(0,0){$3$}}
\put(3,4.5){\makebox(0,0){$4$}}
\put(1,4.5){\makebox(0,0){$5$}}
\end{picture}
\end{minipage}
}
\setlength{\unitlength}{0.8cm}
\subfigure{%
\label{fig:cvbraid}%
\begin{minipage}[t]{0.4\textwidth}
\begin{picture}(3,2.5)(-0.8,-0.4)
\xygraph{!{0;<5mm,0mm>:<0mm,16mm>::}
[]*DL(2){\b{1}}*\cir<3pt>{}
(-[dllll]*U(3){C_{124}}*-{\blacktriangle}
,-[dll]*U(3){C_{135}}*-{\blacktriangle}
,-[d]*U(3){C_{236}}*-{\blacktriangle}
,-[drr]*U(3){C_{456}}*-{\blacktriangle}  
,-[drrrr]*U(3){\Pi}*-{\blacklozenge})
}
\end{picture}
\end{minipage}
}
\caption{\textsf{The braid arrangement $\B$}}
\label{fig:braid}
\end{figure}

\end{example}


\begin{example}[$\operatorname{X}_2$ arrangement] 
\label{ex:x2} 

This is a ``parallel" arrangement (first considered 
by Kohno), for which the LCS formula does not hold, 
see \cite{PY}.

\begin{figure}[ht]
\subfigure{%
\label{fig:dx2}%
\begin{minipage}[t]{0.3\textwidth}
\setlength{\unitlength}{8.5pt}
\begin{picture}(4.5,7.2)(-2,-1)
\multiput(-2,1)(0,2){2}{\line(1,0){8}}
\multiput(1,-2)(2,0){2}{\line(0,1){8}}
\put(-2,4){\line(1,-1){6}}
\put(0,6){\line(1,-1){6}}
\put(4.6,-2.5){\makebox(0,0){$1$}}
\put(6.5,-0.8){\makebox(0,0){$2$}}
\put(6.7,1){\makebox(0,0){$3$}}
\put(6.7,3){\makebox(0,0){$4$}}
\put(3,6.9){\makebox(0,0){$5$}}
\put(1,6.9){\makebox(0,0){$6$}}
\end{picture}
\end{minipage}
}
\setlength{\unitlength}{0.8cm}
\subfigure{%
\label{fig:cvx2}%
\begin{minipage}[t]{0.4\textwidth}
\begin{picture}(3,2.7)(-2,-0.4)
\xygraph{!{0;<9mm,0mm>:<0mm,16mm>::}
[]*DL(2){\b{1}}*\cir<3pt>{}
(-[dll]*U(3){C_{136}}*-{\blacktriangle}
,-[dl]*U(3){C_{245}}*-{\blacktriangle}
,-[d]*U(3){C_{127}}*-{\blacktriangle}
,-[dr]*U(3){C_{347}}*-{\blacktriangle}
,-[drr]*U(3){C_{567}}*-{\blacktriangle}
}
\end{picture}
\end{minipage}
}
\caption{\textsf{The arrangement $\operatorname{X}_2$}}
\label{fig:x2}
\end{figure}

\begin{laundry}

\item 
$Q=xyz(x+y)(x-z)(y-z)(x+y-2z)$.

\item 
$n=7$, $s=11$, $m_2=6$, $m_3=5$.

\item 
$P(X,t)=(1+t)(1+6t+10t^2)$.

\item 
$G^*=G(A_{23},A_{245},A_{35}^{A_{45}},A_{15}^{A_{25}A_{35}A_{45}},
A_{26},A_{46},A_{136}^{A_{23}A_{26}A_{46}},A_{14}^{A_{24}})$. 

\item 
$V_1(G,\K)=C_{136}\cup C_{245}\cup C_{127}\cup C_{347}\cup C_{567}$, 
$V_2(G,\K)=\{ \mathbf{1} \}$. 

\item 
$\beta_{p,1}^{(q)}=\nu_{p,1}=5(p+1)$.

\item 
$b_1(X_N)=5N^{2}+2$.

\item   
$b_1(M_N)=5(N-1)(N-2)$.

\item
$c_1^2(M_N)=N^4(11N^2-36N+29)$, 
$c_2(M_N)=N^4(5N^2-18N+21)$. 

\item 
$\delta_{S_3}=15$, 
$\delta_{A_4}=20$. 

\item 
$a_2=127$, 
$a_3^{\nor}=1,093$,
$a_3=1,138$. 

\item 
$\phi_{1}=7$, $\phi_2=5$, $\phi_3=10$, $\phi_4=15$, $\phi_5=30$, $\phi_6=45$,  
$\phi_7=90$, $\phi_8=150$.  
Conjecture: $\phi_{k}=5 w_k(2)$. 

\item 
$\theta_{1}=7$, $\theta_2=5$, $\theta_3=10$, $\theta_4=15$, 
$\theta_5=20$, $\theta_{k}=5(k-1)$. 

\end{laundry}
\end{example}


\begin{example}[Non-Fano plane]
\label{ex:diamond}

This is a realization of the celebrated non-Fano matroid.  
It is the smallest arrangement for which one of the characteristic 
varieties contains a component not passing through the origin 
(the isolated point $\rho\in V_2(G,\C)$, first found in \cite{CScv}), 
or for which $\RR_{d}(G,\F_q)$ does not coincide with $\RR_{d}(\A)$ mod~$q$, 
for some prime $q$ (as noted in \cite{MS2}).

As we see here, the non-Fano plane is the smallest 
arrangement for which the congruence covers exhibit periodicity. 
Moreover, the LCS ranks follow a particularly 
complicated pattern (the first three were computed by Falk 
and Randell \cite{FR1}, the next three by Peeva \cite{Pe}, 
and the seventh is new).  

\begin{figure}[ht]
\subfigure[Decone]{%
\begin{minipage}[t]{0.25\textwidth}
\setlength{\unitlength}{13pt}
\begin{picture}(6,4.7)(-1.5,-0.7)
\multiput(0,1)(0,2){2}{\line(1,0){4}}
\multiput(1,0)(2,0){2}{\line(0,1){4}}
\put(0,4){\line(1,-1){4}}
\put(0,0){\line(1,1){4}}
\put(4.5,0){\makebox(0,0){$1$}}
\put(4.5,1){\makebox(0,0){$2$}}
\put(4.5,3){\makebox(0,0){$3$}}
\put(4.5,4.5){\makebox(0,0){$4$}}
\put(3,4.5){\makebox(0,0){$5$}}
\put(1,4.5){\makebox(0,0){$6$}}
\end{picture}
\end{minipage}
}
\setlength{\unitlength}{0.8cm}
\subfigure[$V_*(G,\K)$]{%
\label{fig:cvnonfano}%
\begin{minipage}[t]{0.32\textwidth}
\begin{picture}(6,2.2)(-0.5,-1)
\xygraph{!{0;<4.8mm,0mm>:<0mm,16mm>::}
[]*DL(2){\b{1}}*\cir<3pt>{}
(-[dllll]*-{\blacktriangle}
,-[dlll]*-{\blacktriangle}
,-[dll]*-{\blacktriangle}
,-[dl]*-{\blacktriangle}
,-[d] *-{\blacktriangle}
,-[dr]*-{\blacktriangle}  
,-[drr]*U(2.8){\Pi_1}*-{\blacklozenge}
,-[drrr]*U(2.8){\Pi_2}*-{\blacklozenge}
,-[drrrr]*U(2.8){\Pi_3}*-{\blacklozenge}
,[rrr]*DL(2){\rho}*\cir<2pt>{}
(-[dl]
,-[d]
,-[dr])
}
\end{picture}
\end{minipage}
}
\subfigure[$R_*(G,\F_2)$]{%
\label{fig:cvnonfanoz2}%
\begin{minipage}[t]{0.33\textwidth}
\begin{picture}(6,3.2)(-1,-2)
\xygraph{!{0;<5mm,0mm>:<0mm,8mm>::}
[]*DL(2){\b{0}}*\cir<3pt>{}
(-[ddlll]*-{\blacktriangle}
,-[ddll]*-{\blacktriangle}
,-[ddl]*-{\blacktriangle}
,-[dd]*-{\blacktriangle}
,-[ddr] *-{\blacktriangle}
,-[ddrr]*-{\blacktriangle} 
,-[drrr]*L(2.8){\Upsilon'}*{\bigtriangledown}
(
,-[dd]*L(2.8){\Upsilon}*-{\blacksquare}
)
}
\end{picture}
\end{minipage}
}
\caption{\textsf{The non-Fano arrangement}}
\label{fig:nonfano}
\end{figure}

\begin{laundry}

\item 
$Q=xyz(x-y)(x-z)(y-z)(x+y-z)$.

\item 
$n=7$, $s=9$, $m_2=3$, $m_3=6$.

\item 
$P(X,t)=(1+t)(1+3t)^2$.

\item 
$G^*=G(A_{345},  A_{125},  
A_{14}^{A_{34}}, A_{136}, A_{246}^{A_{34}A_{36}})$. 

\item 
$V_1(G,\K)=C_{125}\cup C_{136}\cup C_{246} \cup C_{345} \cup C_{237}\cup C_{567} 
\cup \Pi_1 \cup \Pi_2 \cup \Pi_3$, 
where 
$\Pi_1=C_{(25|36|47)}$, 
$\Pi_2=C_{(17|26|35)}$, 
$\Pi_3=C_{(14|23|56)}$ 
correspond to braid sub-arrangements.
\\[2pt]
$V_{2}(G,\K)=\Pi_1\cap \Pi_2\cap \Pi_3=\{\b{1},\rho \}$, where 
$\rho=(1,-1,-1,1,-1,-1,1)$, unless 
$\ch\K=2$, in which case $V_{2}(G,\K)=\{\b{1}\}$. 

\item 
$\RR_1(G,\F_2)=L_{125}\cup L_{136}\cup L_{246} \cup L_{345} \cup L_{237}\cup L_{567} 
\cup \Upsilon$, where 
\[
\qquad\quad\Upsilon=L_{(2|3|5|6|147)}=
\{(\mu+\nu,\l+\nu,\l+\mu+\nu,\nu,\l+\mu,\l,\mu)\mid 
\l,\mu,\nu\in \F_2\}.
\]
$\RR_2(G,\F_2)=\Upsilon'$, where 
$\Upsilon'=L_{(2|3|5|6)}=\{(0,\l,\l,0,\l,\l,0)\mid \l\in \F_2\}$.
\\[2pt]
$\RR_3(G,\F_2)=\{\b{0}\}$. 

\item 

$\beta_{p,1}^{(q)}=\nu_{p,1}=9(p+1)$, except for:  
$\beta_{2,1}^{(q)}=\nu_{2,1}=24$,  $\beta_{2,2}^{(q)}=\nu_{2,2}=1$.

\item 
$b_1(X_N)=
\begin{cases}
9N^2-3 &\text{if $N$ even,}  \\
9N^2-2 &\text{if $N$ odd.}  
\end{cases}$

\item
$b_1(M_N)=9(N-1)(N-2)$.

\item
$c_1^2(M_N)=N^4(10N^2 - 32N + 25)$,  $c_2(M_N)=N^4(4N^2 - 16N + 21)$.

\item 
$\delta_{S_3}=28$, 
$\delta_{A_4}=36$.   

\item 
$a_2=127$, 
$a_3^{\nor}=1,093$,
$a_3=1,177$. 

\item 
$\phi_{1}=7$, $\phi_2=6$, $\phi_3=17$, $\phi_4=42$,  $\phi_5=123$, $\phi_6=341$, 
$\phi_7=1,041$.

\item 
$\theta_{1}=7$, $\theta_2=6$, $\theta_3=17$, $\theta_4=27$,  $\theta_{k}=9(k-1)$. 

\end{laundry}

\end{example}


\begin{example} [Deleted $\operatorname{B}_3$-arrangement]
\label{ex:deletedB3}

This arrangement is obtained by deleting a plane from 
the reflection arrangement of type  ${\rm B}_3$ (Example \ref{ex:B3}).  
It is a fiber-type arrangement, with exponents $\{1,3,4\}$.  

The deleted $\operatorname{B}_3$-arrangement is 
the smallest arrangement for which one of the characteristic 
varieties contains a positive-dimensional component which does 
{\it not} pass through the origin (see \cite{Su}).  As noted 
in \cite{MS3}, this $1$-dimensional translated torus  
gives rise to $2$-torsion in the homology of $3$-fold 
covers of $X$, and adds $1$ to the number of representations 
of $G$ onto $A_4$.  

Moreover, as we see here, the translated 
component in $V_1$  (together with the two isolated points in $V_2$) 
creates mod~$4$ periodicity in the sequence of first Betti numbers of 
Hirzebruch covering surfaces.  

\begin{figure}[ht]
\subfigure{%
\label{fig:delB3-a}%
\begin{minipage}[t]{0.3\textwidth}
\setlength{\unitlength}{0.65cm}
\begin{picture}(5,4)(1,-0.4)
\multiput(1,0)(1,0){2}{\line(1,1){3}}
\multiput(4,0)(1,0){2}{\line(-1,1){3}}
\multiput(2.5,0)(0.5,0){3}{\line(0,1){3}}
\put(4.2,-0.5){\makebox(0,0){$1$}}
\put(5.2,-0.5){\makebox(0,0){$2$}}
\put(5.2,3.5){\makebox(0,0){$3$}}
\put(4.2,3.5){\makebox(0,0){$4$}}
\put(3.5,3.5){\makebox(0,0){$5$}}
\put(3,3.5){\makebox(0,0){$6$}}
\put(2.5,3.5){\makebox(0,0){$7$}}
\end{picture}
\end{minipage}
}
\subfigure{%
\label{fig:delB3-b}%
\begin{minipage}[t]{0.6\textwidth}
\setlength{\unitlength}{0.6cm}
\begin{picture}(12,4)(0,-3.2)
\xygraph{!{0;<5.5mm,0mm>:<0mm,6.5mm>::}
[]*DL(2){\b{1}}*\cir<3pt>{}
(-[ddddllllll]*{\square}
,-[dddllll]*-{\blacktriangle}
,-[dddlll]*-{\blacktriangle}
,-[dddll]*-{\blacktriangle}
,-[dddl]*-{\blacktriangle}
,-[ddd] *-{\blacktriangle}
,-[dddr]*-{\blacktriangle}  
,-[dddrr]*U(2.5){\Pi_1}*-{\blacklozenge}
,-[dddrrr]*U(2.5){\Pi_2}*-{\blacklozenge}
,-[dddrrrr]*U(2.5){\Pi_3}*-{\blacklozenge}
,-[dddrrrrr]*U(2.5){\Pi_4}*-{\blacklozenge}
,-[dddrrrrrr]*U(2.5){\Pi_5}*-{\blacklozenge}
,[drrrr]*U(2.5){\Omega}*{\bigstar}
,[rrr]*DL(2){\rho_{+}}*\cir<2pt>{}
(-[dddl]
,-[ddd]
,-[dddr]
,-[dr])
,[rrrrr]*DL(2){\rho_{-}}*\cir<2pt>{}  
(-[dddl]
,-[ddd]
,-[dddr]
,-[dl])
)
}
\end{picture}
\end{minipage}
}
\caption{\textsf{The deleted $\operatorname{B}_3$-arrangement}}
\label{fig:delB3arr}
\end{figure}

\begin{laundry}

\item 
$Q=xyz(x-y)(x-z)(y-z)(x-y-z)(x-y+z)$.

\item 
$n=8$, $s=11$, $m_2=4$, $m_3=6$, $m_4=1$.

\item 
$P(X,t)=(1+t)(1+3t)(1+4t)$.

\item 
$G^*=\F_4\rtimes_{\a} \F_3$, where $\a=\{A_{23}, 
A_{13}^{A_{23}}A_{24}, A_{14}^{A_{24}}\}$. 

\item 
$V_1(G,\K)=C_{136}\cup C_{147}\cup C_{235}\cup C_{246}\cup C_{128}\cup 
C_{348}\cup C_{5678}\cup \bigcup_{i=1}^{5} \Pi_i \cup \Omega$, 
where 
$\Omega=\{(t,-t^{-1},-t^{-1},t,t^2,-1,t^{-2},-1) \mid t\in \K^*\}$. 
\\[2pt]
$V_2(G,\K)=C_{5678}\cup \{\rho_{+},\rho_{-}\}$, where 
$\rho_{\pm}=(\pm 1,\mp 1,\mp 1, \pm 1, 1, -1, 1, -1)$, unless 
$\ch\K=2$, in which case $V_{2}(G,\K)=\{\b{1}\}$.  

\item 
$\nu_{p,1}=11(p+1)$, 
$\nu_{p,2}=p^2+p+1$.  \\[2pt]
$\beta_{p,d}^{(q)}=\nu_{p,d}$, 
except for: 
$\beta_{2,1}^{(q)}=27$, 
$\beta_{2,2}^{(q)}=9$,   
and $\beta_{3,1}^{(2)}=45$.  

\item 
$b_1(X_N)=
\begin{cases}
2N^3+11N^2+N-9 &\text{if $N$ even,}  \\
2N^3+11N^2-5   &\text{if $N$ odd,}  
\end{cases}$ 

\item
$b_1(M_N)=
\begin{cases}
(N-1)(2N^2+9N-24) + N-2 &\text{if $N\equiv 0 \mod 4$,}  \\
(N-1)(2N^2+9N-24) + \tfrac{1}{2}(N-2)  &\text{if $N\equiv 2 \mod 4$,}  \\
(N-1)(2N^2+9N-24)  &\text{if $N$ odd.}  
\end{cases}$

\item
$c_1^2(M_N)=N^5(15N^2 - 44N + 31 )$, $c_2(M_N)=2N^5(3N^2- 11N + 13)$. 

\item 
$\delta_{S_3}=63$, 
$\delta_{A_4}=110$.  

\item 
$a_2=255$, 
$a_3^{\nor}=3,280$,
$a_3=3,469$. 

\item 
$\phi_{1}=8$, $\phi_2=9$, $\phi_3=28$, $\phi_4=78$,  
$\phi_k=w_k(3)+w_k(4)$. 

\item 
$\theta_{1}=8$, $\theta_2=9$, $\theta_3=28$, 
$\theta_4=48$,  $\theta_{k}=(k+12)(k-1)$. 

\end{laundry}

\end{example}


\begin{example}[MacLane arrangement]
\label{ex:maclane}

There are two complex realizations of the 
MacLane matroid $\mathtt{ML}_8$ (introduced in \cite{ML}),  
one for each root of the equation $\omega^2+\omega+1=0$.   
The complements of the two arrangements 
are diffeomorphic, but not by an orientation-preserving 
diffeomorphism.  

A presentation for the fundamental group $G$ was 
computed by Rybnikov \cite{Ryb}; the one given here is from \cite{CSbm}. 
The non-local component in $R_2(G,\F_3)$ was first identified in \cite{MS2}. 
Noteworthy is the presence of torsion in the LCS quotients of $G$.

\begin{figure}[ht]
\subfigure[Matroid]{%
\label{fig:matmaclane}%
\begin{minipage}[t]{0.32\textwidth}
\setlength{\unitlength}{0.6cm}
\begin{picture}(3,4.6)(-1.2,-0.5)
\multiput(0,0)(2,0){2}{\line(0,1){4}}
\multiput(0,0)(0,4){2}{\line(1,0){2}}
\multiput(0,0)(0,2){2}{\line(1,1){2}}
\multiput(0,4)(0,-2){2}{\line(1,-1){2}}
\multiput(0,0)(0,2){3}{\circle*{0.3}}
\multiput(1,1)(0,2){2}{\circle*{0.3}}
\multiput(2,0)(0,2){3}{\circle*{0.3}}
\qbezier(1,1)(7,3)(2,4)
\qbezier(1,3)(7,1)(2,0)
\put(0,4){\makebox(0,1){{\small $5$}}}
\put(2,4){\makebox(0,1){{\small $8$}}}
\put(1,1){\makebox(0,1){{\small $2$}}}
\put(0,2){\makebox(-0.8,0){{\small $7$}}}
\put(2,2){\makebox(0.8,0){{\small $6$}}}
\put(1,3){\makebox(0,1){{\small $4$}}}
\put(0,0){\makebox(0,-1){{\small $1$}}}
\put(2,0){\makebox(0,-1){{\small $3$}}}
\end{picture}
\end{minipage}
}
\setlength{\unitlength}{0.5cm}
\subfigure[$V_*(G,\K)$]{%
\label{fig:cvmaclane}%
\begin{minipage}[t]{0.28\textwidth}
\begin{picture}(3,2.3)(-0.5,-1)
\xygraph{!{0;<4mm,0mm>:<0mm,18mm>::}
[]*DL(2){\b{1}}*\cir<3pt>{}
(-[dllll]*U(3){}*-{\blacktriangle}
,-[dlll]*U(3){}*-{\blacktriangle}
,-[dll]*U(3){}*-{\blacktriangle}
,-[dl]*U(3){}*-{\blacktriangle}  
,-[d]*U(3){}*-{\blacktriangle}
,-[dr]*U(3){}*-{\blacktriangle}  
,-[drr]*U(3){}*-{\blacktriangle} 
,-[drrr]*U(3){}*-{\blacktriangle}  
}
\end{picture}
\end{minipage}
}
\subfigure[$R_*(G,\F_3)$]{%
\label{fig:rvmaclane}%
\begin{minipage}[t]{0.3\textwidth}
\begin{picture}(3,2.3)(-0.3,-1)
\xygraph{!{0;<4mm,0mm>:<0mm,18mm>::}
[]*DL(2){\b{0}}*\cir<3pt>{}
(-[dllll]*U(3){}*-{\blacktriangle}
,-[dlll]*U(3){}*-{\blacktriangle}
,-[dll]*U(3){}*-{\blacktriangle}
,-[dl]*U(3){}*-{\blacktriangle}  
,-[d]*U(3){}*-{\blacktriangle}
,-[dr]*U(3){}*-{\blacktriangle}  
,-[drr]*U(3){}*-{\blacktriangle} 
,-[drrr]*U(3){}*-{\blacktriangle}  
,-[drrrr]*U(3.3){\Xi}*-{\blacktriangledown})
}
\end{picture}
\end{minipage}
}
\caption{\textsf{The MacLane arrangement}}
\label{fig:maclane}
\end{figure}

\begin{laundry}

\item 
$Q=xyz(y-x)(z-x)(z+\omega y)(z+\omega^2 x+\omega y)(z-x-\omega^2 y)$.  

\item 
$n=8$, $s=12$, $m_2=4$, $m_3=8$.

\item 
$P(X,t)=(1+t)(1+7t+13t^2)$.

\item 
$G^*=G(A_{456},A_{36},A_{126},A_{134},A_{25},
A_{47}^{A_{57}},A_{157},A_{237}^{A_{57}})$.

\item 
$V_1(G,\K)=C_{126}\cup C_{134}\cup C_{157}\cup C_{237}\cup C_{258}\cup 
C_{368}\cup C_{456}\cup C_{478}$, \\[2pt]
$V_2(G,\K)=\{\b{1}\}$.  

\item 
$\RR_1(G,\F_3)=L_{126}\cup L_{134}\cup L_{157}\cup L_{237}\cup L_{258}\cup 
L_{368}\cup L_{456}\cup L_{478}\cup \Xi$, 
where 
\[
\qquad\quad\Xi=L_{(18|24|35|67)}=
\{(\l-\mu,\l,\mu,-\l,-\mu,\l+\mu,-\l-\mu,\mu-\l)\mid \l,\mu\in\F_3\}.
\]

\item 
$\beta_{p,1}^{(q)}=8(p+1)$ and 
$\nu_{p,1}=8(p+1)$, except for   
$\nu_{3,1}=36$. 

\item 
$b_1(X_N)=8N^2$. 

\item 
$b_1(M_N)=8(N-1)(N-2)$.

\item 
$c_1^2(M_N)=N^5 (  17N^2 - 48N + 32)  $, 
$c_2(M_N)=N^5 ( 7N^2 - 24N +  28 ) $. 

\item 
$\delta_{S_3}=24$, 
$\delta_{A_4}=32$. 

\item 
$a_2=255$, 
$a_3^{\nor}=3,280$,
$a_3=3,352$. 

\item 
$\phi_{1}=8$, $\phi_2=8$, $\phi_3=21$, 
$\phi_4=42$,  $\phi_5=87$,  $\phi_6=105$. 

\item 
$\theta_{1}=8$, $\theta_2=8$, $\theta_3=21$, 
$\theta_4=24$,  $\theta_{k}=8(k-1)$. 

\end{laundry}

\end{example}

\begin{note}
Not all the homology groups of $X_N$ are torsion-free.   
For example, 
$H_1(X_2)=\Z^{32}\oplus \Z_2^4\oplus\Z_4$ and 
$H_1(X_3)=\Z^{72}\oplus \Z_3^8$. 
\end{note}

\begin{note}
Not all the LCS quotients of $G$ are torsion-free.  
For example, $\gr_5 G = \Z^{87} \oplus \Z_2^4\oplus \Z_3$ and 
$\gr_6 G = \Z^{105} \oplus \Z_2^{10}\oplus \Z_4^{21}\oplus 
\Z_{16}^{3}\oplus \Z_3^{16}$.
\end{note}


\begin{example}[$\operatorname{B}_3$-arrangement]
\label{ex:B3}

This is the complexification of the reflection arrangement 
of type $\operatorname{B}_3$, consisting of the nine planes 
of symmetry of the cube in $\R^3$, with vertices at $(\pm 1,\pm 1, \pm 1)$.  
It is a fiber-type arrangement, with exponents $\{1,3,5\}$.  
The resonance variety has an essential component,   
first identified in \cite{Fa}.  

\begin{figure}[ht]
\subfigure{%
\label{fig:B3-a}%
\begin{minipage}[t]{0.3\textwidth}
\setlength{\unitlength}{0.58cm}
\begin{picture}(5,4.8)(-0.2,-1)
\multiput(1,0)(1,0){2}{\line(1,1){3}}
\multiput(4,0)(1,0){2}{\line(-1,1){3}}
\multiput(2.5,0)(0.5,0){3}{\line(0,1){3}}
\put(1,1.5){\line(1,0){4}}
\put(4.2,-0.5){\makebox(0,0){$1$}}
\put(5.2,-0.5){\makebox(0,0){$2$}}
\put(5.5,1.5){\makebox(0,0){$3$}}
\put(5.2,3.5){\makebox(0,0){$4$}}
\put(4.2,3.5){\makebox(0,0){$5$}}
\put(3.5,3.5){\makebox(0,0){$6$}}
\put(3,3.5){\makebox(0,0){$7$}}
\put(2.5,3.5){\makebox(0,0){$8$}}
\end{picture}
\end{minipage}
}
\subfigure{%
\label{fig:B3-b}%
\begin{minipage}[t]{0.6\textwidth}
\setlength{\unitlength}{0.6cm}
\begin{picture}(6,4.8)(0.5,-3.5)
\xygraph{!{0;<3.8mm,0mm>:<0mm,6.5mm>::}
[]*DL(2){\b{1}}*\cir<3pt>{}
(-[ddddllllllllll]*{\square}
,-[ddddlllllllll]*{\square}
,-[ddddllllllll]*{\square}
,-[dddlllll]*-{\blacktriangle}
,-[dddllll]*-{\blacktriangle}
,-[dddlll]*-{\blacktriangle}
,-[dddll]*-{\blacktriangle} 
,-[dddl]*-{\blacklozenge}
,-[ddd]*-{\blacklozenge}
,-[dddr]*-{\blacklozenge}
,-[dddrr]*-{\blacklozenge}
,-[dddrrr]*-{\blacklozenge}
,-[dddrrrr]*-{\blacklozenge}
,-[dddrrrrr]*-{\blacklozenge}
,-[dddrrrrrr]*-{\blacklozenge}
,-[dddrrrrrrr]*-{\blacklozenge}
,-[dddrrrrrrrr]*-{\blacklozenge}
,-[dddrrrrrrrrr]*-{\blacklozenge}
,-[dddrrrrrrrrrr]*U(2.9){\Gamma}*-{\clubsuit}
,[rrrr]*DL(1.8){\sml{\rho_{+}}}*\cir<2pt>{}
(-[dddl]
,-[ddd]
,-[dddrr]
,-[dddrrrrrr])
,[rrrrrr]*DL(1.8){\sml{\rho_{-}}}*\cir<2pt>{}  
(-[dddl]
,-[ddd]
,-[dddr]
,-[dddrrrr])
,[rrrrrrrr]*DL(1.6){\sml{\rho_{+}\rho_{-}}}*\cir<2pt>{}  
(-[dddll]
,-[ddd]
,-[dddr]
,-[dddrr])
)
}
\end{picture}
\end{minipage}
}
\caption{\textsf{The $\operatorname{B}_3$-arrangement}}
\label{fig:B3arr}
\end{figure}

\begin{laundry}

\item 
$Q=xyz(x-y)(x-z)(y-z)(x-y-z)(x-y+z)(x+y-z)$.

\item 
$n=9$, $s=13$, $m_2=6$, $m_3=4$, $m_4=3$.

\item 
$P(X,t)=(1+t)(1+3t)(1+5t)$.

\item 
$G^*=\F_5\rtimes_{\a} \F_3$, where $\a=\{A_{234}, 
A_{14}^{A_{24}A_{34}}A_{25},A_{35}^{A_{23}A_{25}}\}$. 

\item 

$V_1(G,\K)=C_{147}\cup C_{257}\cup C_{129}\cup C_{459}\cup 
C_{2346}\cup C_{1358}\cup C_{6789}\cup \bigcup_{i=1}^{11}\Pi_i 
\cup \Gamma$,~where
\[
\qquad\quad\Gamma=C_{(156 | 248 | 379)}=
\{ (t,s,(st)^{-2},s,t,t^2,(st)^{-1},s^2,(st)^{-1})
\mid s,t\in \K^*\}.
\]
$V_2(G,\K)=C_{2346}\cup C_{1358}\cup C_{6789}\cup 
\{\rho_{+},\rho_{-}, \rho_{+}\rho_{-}\}$, where 
\[
\rho_{\pm}=(\pm 1,\mp 1,1, \mp 1,  \pm 1,1,-1,1,-1).
\]

\item 
$\nu_{p,1}=16(p+1)$,  
$\nu_{p,2}=3(p^2+p+1)$.  \\[2pt] 
$\beta_{p,d}^{(q)}=\nu_{p,d}$, 
except for:  
$\beta_{2,1}^{(q)}=36$,  
$\beta_{2,2}^{(q)}=24$.  

\item 
$b_1(X_N)=
\begin{cases}
6N^3+16N^2-19 &\text{if $N$ even,}  \\
6N^3+16N^2-13 &\text{if $N$ odd.}  
\end{cases}$

\item
$b_1(M_N)=2(N-1)(3N^2+5N-19)$.

\item
$c_1^2(M_N)=2 N^6 (10N^2-28N+19)$, 
$c_2(M_N)=2 N^6 (4N^2-14 N+15)$.  

\item 
$\delta_{S_3}=132$, 
$\delta_{A_4}=259$. 

\item 
$a_2=511$, 
$a_3^{\nor}=9,841$, 
$a_3=10,237$. 

\item 
$\phi_{1}=9$, $\phi_2=13$, $\phi_3=48$, $\phi_4=168$, 
$\phi_k=w_k(3)+w_k(5)$. 

\item 
$\theta_{1}=9$, $\theta_2=13$, $\theta_3=48$, 
$\theta_4=93$,  $\theta_{k}=(3k+19)(k-1)$. 

\end{laundry}

\end{example}


\begin{example}[Pappus arrangement] 
\label{ex:pappus}

This is a realization of the classical $(9_3)_1$  
configuration of Hilbert and Cohn-Vossen~\cite{HC}.   

\begin{figure}[ht]
\subfigure{%
\label{fig:dpap}%
\begin{minipage}[t]{0.3\textwidth}
\setlength{\unitlength}{8.5pt}
\begin{picture}(5,9)(-2,-1)
\multiput(-1.5,1)(0,2){2}{\line(1,0){9}}
\multiput(1,-1)(2,0){2}{\line(0,1){10}}
\put(-1,7){\line(3,-4){6}}
\put(0,9){\line(3,-4){7}}
\put(0,-1){\line(1,2){4.5}}
\put(-1.5,6){\line(3,-2){9}}
\put(-0.35,-1.8){\makebox(0,0){$1$}}
\put(1,-1.8){\makebox(0,0){$2$}}
\put(3,-1.8){\makebox(0,0){$3$}}
\put(5.15,-1.8){\makebox(0,0){$4$}}
\put(7.2,-1.1){\makebox(0,0){$5$}}
\put(8.05,-0.4){\makebox(0,0){$6$}}
\put(8.15,1){\makebox(0,0){$7$}}
\put(8.15,3){\makebox(0,0){$8$}}
\end{picture}
\end{minipage}
}
\setlength{\unitlength}{0.8cm}
\subfigure{%
\label{fig:cvpap}%
\begin{minipage}[t]{0.4\textwidth}
\begin{picture}(3,2.5)(-1,-0.4)
\xygraph{!{0;<5.5mm,0mm>:<0mm,18mm>::}
[]*DL(2){\b{1}}*\cir<3pt>{}
(-[dllll]*-{\blacktriangle}
,-[dlll]*-{\blacktriangle}
,-[dll]*-{\blacktriangle}
,-[dl]*-{\blacktriangle}
,-[d]*-{\blacktriangle}
,-[dr]*-{\blacktriangle}
,-[drr]*-{\blacktriangle}
,-[drrr]*-{\blacktriangle}
,-[drrrr]*-{\blacktriangle}
,-[drrrrr]*U(2.9){\Psi}*-{\spadesuit}
}
\end{picture}
\end{minipage}
}
\caption{\textsf{The Pappus arrangement $(9_3)_1$}}
\label{fig:pap}
\end{figure}

\begin{laundry}

\item 
$Q=xyz(x-y)(y-z)(x-y-z)(2x+y+z)(2x+y-z)(2x-5y+z)$. 

\item 
$n=9$, $s=18$, $m_2=9$, $m_3=9$.

\item 
$P(X,t)=(1+t)(1+8t+19t^2)$.

\item 
$G^*=G(A_{567},A_{47},A_{37},A_{127},A_{34},A_{58},A_{368}^{A_{56}},
A_{148}^{A_{34}A_{38}A_{58}A_{68}},A_{16}^{A_{36}A_{56}}$,\\ $
A_{28}^{A_{38}A_{48}A_{58}A_{68}},A_{135},A_{246}^{A_{34}A_{36}A_{56}},
A_{25}^{A_{35}})$.

\item 
$V_1(G,\K)=C_{567}\cup C_{127}\cup C_{368}\cup C_{148}
\cup C_{135}\cup C_{246}\cup 
C_{239} \cup C_{459}\cup C_{789}\cup \Psi$, where
$\Psi=C_{(169 | 258 | 347)}=
\{ (t,s,(st)^{-1},(st)^{-1},s,t,(st)^{-1},s,t)
\mid s,t\in \K^*\}.
$
$V_{2}(G,\K)=\{\b{1}\}$.

\item 
$\beta_{p,1}^{(q)}=\nu_{p,1}=10(p+1)$. 

\item 
$b_1(X_N)=10N^2-1$. 

\item
$b_1(M_N)=10(N-1)(N-2)$.

\item 
$c_1^2(M_N)=9 N^6 (N-1) (3N-5)$, 
$c_2(M_N)= 12 N^6 (N^2-3N+3)$. 

\item 
$\delta_{S_3}=30$, 
$\delta_{A_4}=40$. 

\item 
$a_2=511$, 
$a_3^{\nor}=9,841$,
$a_3=9,931$. 

\item 
$\phi_{1}=9$, $\phi_2=9$, $\phi_3=20$, 
$\phi_4=30$,  $\phi_5=60$, $\phi_6=90$, $\phi_7=180$, $\phi_8=300$. \\
Conjecture:  $\phi_k=10w_k(2)$.  

\item 
$\theta_{1}=9$, $\theta_2=9$, $\theta_3=20$, 
$\theta_4=30$, $\theta_5=50$, $\theta_{k}=10(k-1)$ 

\end{laundry}

\end{example}


\begin{example}[Non-Pappus arrangement] 
\label{ex:nonpappus}

This is a realization of the $(9_3)_2$ configuration from \cite{HC}.  
It has the same number of multiple points as the Pappus arrangement, 
but the difference in their position is reflected in several  
invariants of the complement. 

\begin{figure}[ht]
\subfigure{%
\label{fig:dnonpap}%
\begin{minipage}[t]{0.3\textwidth}
\setlength{\unitlength}{9pt}
\begin{picture}(5,7.7)(-2,-1.35)
\multiput(-1,1.5)(0,1){2}{\line(1,0){9.2}}
\multiput(1,-1.5)(3,0){2}{\line(0,1){8}}
\put(-1,3.5){\line(2,-1){9.2}}
\put(-1,4.5){\line(2,-1){9.2}}
\put(-0.95,4.8){\line(3,-2){9.2}}
\put(0,6.5){\line(1,-1){8}}
\put(1,-2.2){\makebox(0,0){$1$}}
\put(4,-2.2){\makebox(0,0){$2$}}
\put(8,-2.2){\makebox(0,0){$3$}}
\put(8.85,-2.1){\makebox(0,0){$4$}}
\put(9,-1.2){\makebox(0,0){$5$}}
\put(9,-0.3){\makebox(0,0){$6$}}
\put(9,1.5){\makebox(0,0){$7$}}
\put(9,2.5){\makebox(0,0){$8$}}
\end{picture}
\end{minipage}
}
\setlength{\unitlength}{0.8cm}
\subfigure{%
\label{fig:cvnonpap}%
\begin{minipage}[t]{0.4\textwidth}
\begin{picture}(3,2.3)(-1.3,0)
\xygraph{!{0;<5.5mm,0mm>:<0mm,18mm>::}
[]*DL(2){\b{1}}*\cir<3pt>{}
(-[dllll]*-{\blacktriangle}
,-[dlll]*-{\blacktriangle}
,-[dll]*-{\blacktriangle}
,-[dl]*-{\blacktriangle}
,-[d]*-{\blacktriangle}
,-[dr]*-{\blacktriangle}
,-[drr]*-{\blacktriangle}
,-[drrr]*-{\blacktriangle}
,-[drrrr]*-{\blacktriangle}
}
\end{picture}
\end{minipage}
}
\caption{\textsf{The non-Pappus arrangement $(9_3)_2$}}
\label{fig:nonpap}
\end{figure}

\begin{laundry}

\item 
$Q=xyz(x+y)(y+z)(x+3z)(x+2y+z)(x+2y+3z)(2x+3y+3z)$.

\item 
$n=9$, $s=18$, $m_2=9$, $m_3=9$.

\item 
$P(X,t)=(1+t)(1+8t+19t^2)$.

\item 
$G^*=G(A_{345},A_{367},A_{25},A_{247}^{A_{34}A_{36}^{-1}},A_{26}^{A_{36}},A_{238},
A_{57}^{A_{67}},A_{68},A_{48},A_{17}^{A_{27}A_{37}A_{47}A_{67}}$,\\
$A_{158},A_{146}^{A_{24}A_{34}A_{36}},A_{13}^{A_{23}})$. 

\item 
$V_1(G,\K)=C_{345}\cup C_{367}\cup C_{247}\cup C_{238}\cup C_{158}\cup C_{146}\cup 
C_{129} \cup C_{569}\cup C_{789}$, 
$V_{2}(G,\K)=\{\b{1}\}$.

\item 
$\beta_{p,1}^{(q)}=\nu_{p,1}=9(p+1)$. 

\item 
$b_1(X_N) = 9N^2$. 

\item
$b_1(M_N) = 9(N-1)(N-2)$.

\item
$c_1^2(M_N) =9 N^6 (N-1) (3N-5)$, 
$c_2(M_N) = 12 N^6 (N^2-3N+3)$.  

\item 
$\delta_{S_3}=27$, 
$\delta_{A_4}=36$. 

\item 
$a_2=511$, 
$a_3^{\nor}=9,841$,
$a_3=9,922$. 

\item 
$\phi_{1}=9$, $\phi_2=9$, $\phi_3=18$, 
$\phi_4=27$,  $\phi_5=54$, $\phi_6=81$, 
$\phi_7=162$, $\phi_8=270$. \\
Conjecture:  $\phi_k=9 w_k(2)$.  

\item 
$\theta_{1}=9$, $\theta_2=9$, $\theta_3=18$, 
$\theta_4=27$, $\theta_5=36$, $\theta_{k}=9(k-1)$.

\end{laundry}

\end{example}


\begin{example}  [Ziegler arrangements]
\label{ex:ziegler}

These arrangements, $\A_1$ and $\A_2$, are both fiber-type, 
with exponents $\{1, 6, 6\}$.  They 
are combinatorially very close, yet a subtle difference 
in the respective lattices gets detected by some 
of the more refined topological invariants of 
the complements, such as the number, $\tau$, of translated 
components in the characteristic varieties \cite{Su},   
or the Hall invariant  $\delta_{A_4}$. 
Furthermore, the sequences $b_1(X_N(\A_i))$ and $b_1(M_N(\A_i))$ 
are polynomially periodic (of period $2$ and $4$, respectively), 
but the corresponding polynomials differ for $N$ even.

\begin{figure}[ht]
\subfigure{%
\label{fig:decz1}%
\begin{minipage}[t]{0.24\textwidth}
\setlength{\unitlength}{0.3cm}
\begin{picture}(4,9.2)(0.4,-2.2)
\multiput(2,1)(0,-1){3}{\line(1,1){6}}
\put(2,-3){\line(1,1){6}}
\multiput(2.5,-3)(1,0){6}{\line(0,1){10}}
\multiput(2,1.5)(0,1){2}{\line(1,0){6}}
\put(8.7,1.5){\makebox(0,0){$\sml{1}$}}
\put(8.7,2.5){\makebox(0,0){$\sml{2}$}}
\put(8.7,3.4){\makebox(0,0){$\sml{3}$}}
\put(8.7,5){\makebox(0,0){$\sml{4}$}}
\put(8.7,6){\makebox(0,0){$\sml{5}$}}
\put(8.7,7.2){\makebox(0,0){$\sml{6}$}}
\put(7.5,8){\makebox(0,0){$\sml{7}$}}
\put(6.5,8){\makebox(0,0){$\sml{8}$}}
\put(5.5,8){\makebox(0,0){$\sml{9}$}}
\put(4.5,8){\makebox(0,0){$\sml{10}$}}
\put(3.5,8){\makebox(0,0){$\sml{11}$}}
\put(2.3,8){\makebox(0,0){$\sml{12}$}}
\end{picture}
\end{minipage}
}
\subfigure{%
\label{fig:cvz1}%
\begin{minipage}[t]{0.65\textwidth}
\setlength{\unitlength}{0.7cm}
\begin{picture}(6,4)(0.8,-3.2)
\xygraph{!{0;<2.5mm,0mm>:<0mm,5mm>::}
[]*DL(2){\b{1}}*\cir<3pt>{}
(-[ddddddlllllllllllllllllll]*!{\boxtimes}
,-[dddddllllllllllllll]*!{\boxplus}
,-[ddddllllllllll]*-{\blacktriangle}
,-[ddddlllllllll]*-{\blacktriangle}
,-[ddddllllllll]*-{\blacktriangle}
,-[ddddlllllll]*-{\blacktriangle} 
,-[ddddllllll]*-{\blacktriangle}
,-[ddddlllll]*-{\blacktriangle} 
,-[ddddllll]*-{\blacktriangle} 
,-[ddddlll]*-{\blacktriangle}
,-[ddddll]*-{\blacktriangle} 
,-[ddddl]*-{\blacklozenge}
,-[dddd]*-{\blacklozenge}
,-[ddddr]*-{\blacklozenge}
,-[ddddrr]*-{\blacklozenge}
,-[ddddrrr]*-{\blacklozenge}
,-[ddddrrrr]*-{\blacklozenge}
,-[ddddrrrrr]*-{\blacklozenge}
,-[ddddrrrrrr]*-{\blacklozenge}
,-[ddddrrrrrrr]*-{\blacklozenge}
,-[ddddrrrrrrrr]*-{\blacklozenge}
,-[ddddrrrrrrrrr]*-{\blacklozenge}
,-[ddddrrrrrrrrrr]*-{\blacklozenge}
,-[ddddrrrrrrrrrrr]*-{\blacklozenge}
,-[ddddrrrrrrrrrrrr]*-{\blacklozenge}
,-[ddddrrrrrrrrrrrrr]*-{\blacklozenge}
,-[ddddrrrrrrrrrrrrrr]*-{\blacklozenge}
,-[ddddrrrrrrrrrrrrrrr]*-{\blacklozenge}
,-[ddddrrrrrrrrrrrrrrrr]*-{\blacklozenge}
,[drrrrrrrrrr]*{\bigstar}
,[drrrrrrrrrrrr]*{\bigstar}
,[drrrrrrrrrrrrrr]*{\bigstar}
,[rrrrrrrrr]*\cir<2pt>{}
(-[dr]
,-[ddddl]
,-[dddd]
,-[ddddr])
,[rrrrrrrrrrr]*\cir<2pt>{}  
(-[dl]
,-[dr]
,-[ddddl]
,-[dddd]
,-[ddddr])
,[rrrrrrrrrrrrr]*\cir<2pt>{}  
(-[dl]
,-[dr]
,-[ddddl]
,-[dddd]
,-[ddddr])
,[rrrrrrrrrrrrrrr]*\cir<2pt>{}  
(-[dl]
,-[ddddl]
,-[dddd]
,-[ddddr])
)
}
\end{picture}
\end{minipage}
}
\caption{\textsf{The first Ziegler arrangement}}
\label{fig:ziegler1}
\end{figure}

\begin{figure}[ht]
\subfigure{%
\label{fig:decz2}%
\begin{minipage}[t]{0.24\textwidth}
\setlength{\unitlength}{0.3cm}
\begin{picture}(4,9.2)(0.4,-2.2)
\multiput(2,1)(0,-1){2}{\line(1,1){6}}
\multiput(2,-2)(0,-1){2}{\line(1,1){6}}
\multiput(2.5,-3)(1,0){6}{\line(0,1){10}}
\multiput(2,1.5)(0,1){2}{\line(1,0){6}}
\put(8.7,1.5){\makebox(0,0){$\sml{1}$}}
\put(8.7,2.5){\makebox(0,0){$\sml{2}$}}
\put(8.7,3.4){\makebox(0,0){$\sml{3}$}}
\put(8.7,5){\makebox(0,0){$\sml{4}$}}
\put(8.7,6){\makebox(0,0){$\sml{5}$}}
\put(8.7,7.2){\makebox(0,0){$\sml{6}$}}
\put(7.5,8){\makebox(0,0){$\sml{7}$}}
\put(6.5,8){\makebox(0,0){$\sml{8}$}}
\put(5.5,8){\makebox(0,0){$\sml{9}$}}
\put(4.5,8){\makebox(0,0){$\sml{10}$}}
\put(3.5,8){\makebox(0,0){$\sml{11}$}}
\put(2.3,8){\makebox(0,0){$\sml{12}$}}
\end{picture}
\end{minipage}
}
\subfigure{%
\label{fig:cvz2}%
\begin{minipage}[t]{0.65\textwidth}
\setlength{\unitlength}{0.7cm}
\begin{picture}(6,4)(0.8,-3.2)
\xygraph{!{0;<2.5mm,0mm>:<0mm,5mm>::}
[]*DL(2){\b{1}}*\cir<3pt>{}
(-[ddddddlllllllllllllllllll]*!{\boxtimes}
,-[dddddllllllllllllll]*!{\boxplus}
,-[ddddllllllllll]*-{\blacktriangle}
,-[ddddlllllllll]*-{\blacktriangle}
,-[ddddllllllll]*-{\blacktriangle}
,-[ddddlllllll]*-{\blacktriangle} 
,-[ddddllllll]*-{\blacktriangle}
,-[ddddlllll]*-{\blacktriangle} 
,-[ddddllll]*-{\blacktriangle} 
,-[ddddlll]*-{\blacktriangle}
,-[ddddll]*-{\blacktriangle} 
,-[ddddl]*-{\blacklozenge}
,-[dddd]*-{\blacklozenge}
,-[ddddr]*-{\blacklozenge}
,-[ddddrr]*-{\blacklozenge}
,-[ddddrrr]*-{\blacklozenge}
,-[ddddrrrr]*-{\blacklozenge}
,-[ddddrrrrr]*-{\blacklozenge}
,-[ddddrrrrrr]*-{\blacklozenge}
,-[ddddrrrrrrr]*-{\blacklozenge}
,-[ddddrrrrrrrr]*-{\blacklozenge}
,-[ddddrrrrrrrrr]*-{\blacklozenge}
,-[ddddrrrrrrrrrr]*-{\blacklozenge}
,-[ddddrrrrrrrrrrr]*-{\blacklozenge}
,-[ddddrrrrrrrrrrrr]*-{\blacklozenge}
,-[ddddrrrrrrrrrrrrr]*-{\blacklozenge}
,-[ddddrrrrrrrrrrrrrr]*-{\blacklozenge}
,-[ddddrrrrrrrrrrrrrrr]*-{\blacklozenge}
,-[ddddrrrrrrrrrrrrrrrr]*-{\blacklozenge}
,[drrrrrrrrr]*{\bigstar}
,[drrrrrrrrrrrrrr]*{\bigstar}
,[rrrrrrrr]*\cir<2pt>{}
(-[dr]
,-[ddddl]
,-[dddd]
,-[ddddr])
,[rrrrrrrrrr]*\cir<2pt>{}  
(-[dl]
,-[ddddl]
,-[dddd]
,-[ddddr])
,[rrrrrrrrrrrrr]*\cir<2pt>{}  
(-[dr]
,-[ddddl]
,-[dddd]
,-[ddddr])
,[rrrrrrrrrrrrrrr]*\cir<2pt>{}  
(-[dl]
,-[ddddl]
,-[dddd]
,-[ddddr])
)
}
\end{picture}
\end{minipage}
}
\caption{\textsf{The second Ziegler arrangement}}
\label{fig:ziegler2}
\end{figure}

\begin{laundry}

\item 
$Q_{\A_1}=xyz(x-y)(y-z)(x-z)(x-2z)(x-3z)(x-4z)(x-5z)\times$ \\ 
\mbox{\qquad\qquad} $(x-y-z)(x-y-2z)(x-y-4z)$, \\[3pt]
$Q_{\A_2}=xyz(x-y)(y-z)(x-z)(x-2z)(x-3z)(x-4z)(x-5z)\times$ \\ 
\mbox{\qquad\qquad} $(x-y-z)(x-y-3z)(x-y-4z)$.

\item 
$n=13$, $s=31$, $m_2=20$, $m_3=9$, $m_5=1$, $m_7=1$.

\item 
$P(X_i,t)=(1+t)(1+6t)^2$.

\item 
$G_i^*=F_6\rtimes_{\a_i} F_6$, where 
$\a_1=\{A_{23}, A_{13}^{A_{23}}, A_{24}, A_{14}^{A_{24}}A_{25},  
 A_{15}^{A_{25}}A_{26}, A_{16}^{A_{26}} \}$, 
$\a_2=\{A_{23}, A_{13}^{A_{23}}A_{24}, A_{14}^{A_{24}},  A_{25}, 
 A_{15}^{A_{25}}A_{26}, A_{16}^{A_{26}} \}$. 

\item 
Each variety $V_1(G_i,\K)$ has $11$ local components, 
and $18$ components corresponding to braid sub-arrangements.  
In addition, $V_1(G_1,\K)$ has $\tau_1=3$ components
corresponding to deleted $\operatorname{B}_3$ sub-arrangements,
whereas $V_1(G_2,\K)$ has $\tau_2=2$ such components (not 
passing through the origin). 
\\[3pt]
$V_2(G_i,\K)=\{\b{1},\rho_1,\rho_2,\rho_3,\rho_4\}$, with $\rho_j^2=\b{1}$;\  
$V_2(G_i,\K)=\{\b{1}\}$, if $\ch\K=2$. 

\item 
$\nu_{p,1}=27(p+1)$, $\nu_{p,3}=\tfrac{p^4-1}{p-1}$, 
$\nu_{p,5}=\tfrac{p^6-1}{p-1}$.  \\[3pt]
$\beta_{p,d}^{(q)}=\nu_{p,d}$, except for 
$\beta_{2,1}^{(q)}=69$, 
$\beta_{2,2}^{(q)}=4$, and
\[
\beta_{3,1}^{(2)}(G_1)=111,\quad 
\beta_{3,1}^{(2)}(G_2)=110. 
\]

\item 
$b_1(X_N(\A_i))=
\begin{cases}
5 N^6 + 3 N^4 + 27 N^2 +\tau_i(N-2) - 26 &\text{if $N$ even,}  \\
5 N^6 + 3 N^4 + 27 N^2 -22 &\text{if $N$ odd.}  
\end{cases}$ 

\item
$b_1(M_N(\A_i))=
\begin{cases}
f(N) + \tau_i(N-2) &\text{if $N\equiv 0 \mod 4$,}  \\
f(N) + \tfrac{\tau_i}{2}(N-2)  &\text{if $N\equiv 2 \mod 4$,}  \\
f(N) &\text{if $N$ odd,}  
\end{cases}$\\[3pt]
where $f(N)= (N-1)(5N^5-2N^4+N^3-4N^2+23N-58)$.

\item
$c_1^2(M_N(\A_i))=N^{10} (  57 N^2- 140 N + 81)$,\\[2pt]
$c_2(M_N(\A_i))= N^{10} ( 25 N^2- 70 N + 59)$. 

\item 
$\delta_{S_3}=7,903$,  and 
\[\delta_{A_4}(G_1)=124,435, \quad
\delta_{A_4}(G_2)=124,434.
\] 

\item 
$a_2=8,191$, 
$a_3^{\nor}=797,161$,
$a_3=820,870$. 

\item 
$\phi_{1}=13$, $\phi_2=30$, $\phi_3=140$, 
$\phi_4=630$,  $\phi_k=2w_k(6)$. 

\item 
$\theta_{1}=13$, $\theta_2=30$, $\theta_3=140$, $\theta_4=336$,  
$\theta_{k}=\frac{(k-1)(k^4+10k^3+47k^2+86k+696)}{24}$. 

\end{laundry}

\end{example}


\end{document}